\documentclass[11pt]{amsart}
\usepackage{amssymb}
\usepackage{amsfonts}
\usepackage{amscd}
\usepackage{amsmath}
\usepackage{mathrsfs}
\usepackage{curves}
\usepackage{epsfig}
\usepackage[pass]{geometry}
\usepackage{graphicx}
\usepackage{verbatim}
\usepackage{latexsym}
\usepackage{eucal, tikz}
\usetikzlibrary{matrix}

\newtheorem{theorem}{Theorem}[section]
\newtheorem{lemma}[theorem]{Lemma}
\newtheorem{prop}[theorem]{Proposition}
\newtheorem{cor}[theorem]{Corollary}

\newtheorem{definition}[theorem]{Definition}

\def \mcb {{\mathcal B}}
\def \mcc {{\mathcal C}}
\def \mcd {{\mathcal D}}
\def \mce {{\mathcal E}}
\def \mcf {{\mathcal F}}

\def \mci {{\mathcal I}}
\def \mcj {{\mathcal J}}
\def \mck {{\mathcal K}}
\def \mcl {{\mathcal L}}
\def \mcm {{\mathcal M}}
\def \mcn {{\mathcal N}}

\def \mcr {{\mathcal R}}

\def \mcu {{\mathcal U}}
\def \mcv {{\mathcal V}}
\def \mcx {\mathcal{X}}

\def \mcz {{\mathcal Z}}

\def \mbr {{\mathbb R}}
\def \mbs {{\mathbb S}}

\def \id {\operatorname{Id}}
\def \comp {\operatorname{comp}}
\def \loc {\operatorname{loc}}

\def \re {\operatorname{Re}}
\def \diag{\operatorname{Diag}}
\def \supp {\operatorname{supp }}

\def \beqq {\begin{equation}}
\def \eeqq {\end{equation}}
\def \WF {\operatorname{WF}}

\def \bpf {\begin{proof}}
\def \epf {\end{proof}}
\def \beq {\begin{equation*}}
\def \eeq {\end{equation*}}

\def \La {\Lambda}

\def \p {\partial}
\def \ha {\frac{1}{2}}

\def \tilde {\widetilde}
\def \dist {\operatorname{dist}}
\def \vol {\text{Vol}}

\begin{document}
\title[light ray transform]{Microlocal analysis of the light ray transform on globally hyperbolic Lorentzian manifolds}
\author{Yiran Wang}
\address{Yiran Wang
\newline
\indent Department of Mathematics, Emory University}
\email{yiran.wang@emory.edu}
\begin{abstract}
For the light ray transform on globally hyperbolic Lorentzian manifolds of dimension $n+1 \geq 3$ acting on compactly supported distributions, we show that the Schwartz kernel of the normal operator is a paired Lagrangian distribution  with non-vanishing principal symbols on each Lagrangians. We obtain Sobolev estimates for the light ray transform, and clarify the determination of light-like singularities using the normal operator.  
\end{abstract}
\date{\today}
 
\maketitle

\section{Introduction}
Let $(\mcm, g)$ be an $n+1, n\geq 2$ dimensional smooth Lorentzian manifold. Let $\gamma$ be a complete light-like (or null) geodesic which means that $\gamma(s)$ is defined for $s\in \mbr$ and $\dot \gamma(s)$ satisfies $g(\dot \gamma(s), \dot \gamma(s)) = 0$. We consider the light ray transform 
\beqq\label{eq-lray0}
(Lf)(\gamma) = \int_{\mbr}  f(\gamma(s)) ds, \quad f\in C_0^\infty(\mcm)
\eeqq
when the integral is well-defined. Note that even for $C_0^\infty$ functions, the integral may not converge because $\gamma$ may be trapped in the support of $f$.  
In this work, we study the transform on globally hyperbolic Lorentzian manifold $(\mcm, g)$.  We recall from Bernal and Sanchez \cite[Theorem 1.1]{BeSa} that  such $(\mcm, g)$ can be identified with  
\beqq\label{eq-metric}
 \mbr\times \mcn, \quad g = -\beta(t, x)dt^2 + h(t, x, dx)
\eeqq
where $\mcn$ is a smooth $n$ dimensional spacelike Cauchy hyper-surface, $\beta: \mbr\times \mcn \rightarrow (0, \infty)$ is a smooth function, and $h$ is a smooth family of Riemannian metrics on $\mcn.$  We assume that $(\mcm, g)$ has no conjugate points and is null-geodesically complete. We remark that these global assumptions are imposed to simplify the elaboration. Our results apply when these assumptions hold locally on the support of $f$ as those studied in \cite{LOSU2}.

We introduce some notions to state the main result. It is known that light-like geodesics are invariant under conformal diffeomorphisms of $(\mcm, g)$, however the light ray transform \eqref{eq-lray0} is not invariant.  For globally hyperbolic manifolds, there is a natural parametrization of the light ray transform using the unit sphere bundle  of $(\mcn, g|_\mcn)$ denoted by $\mcc$, see \eqref{eq-lorlray}.  We will specify the measure on $\mcc$ later. We consider $L: C_0^\infty(\mcm)\rightarrow C^\infty(\mcc)$, and its transpose $L^t: C_0^\infty(\mcc)\rightarrow C^\infty(\mcm)$.  Our main result is about the microlocal structure of the normal operator $N = L^tL.$ The notions for paired Lagrangians will be recalled in Section \ref{sec-pair}.  
\begin{theorem}\label{thm-main}
Let $(\mcm, g)$ be a globally hyperbolic Lorentzian manifold of dimension $n+1, n\geq 2$. Suppose $(\mcm, g)$ is null-geodesic complete without conjugate points. Consider the normal operator $N = L^t L$ of the light ray transform $L$. Then the Schwartz kernel $K_N  \in I^{-n/2, n/2- 1}(\mcm\times\mcm; \La_0, \La_1)$, in which $\La_0, \La_1$ are two cleanly intersection Lagrangians defined in \eqref{eq-laggen}. Let $\Sigma = \La_0\cap \La_1.$ The principal symbols  of $K_N$ on $\La_0\backslash \Sigma, \La_1\backslash \Sigma$ are non-vanishing. 
\end{theorem}
As a consequence, we obtain Sobolev estimates for the light ray transform. We remark that for the Minkowski spacetime, related estimates were obtained by Greenleaf and Seeger \cite{GrSe}.  
\begin{theorem}\label{thm-mainest}
Consider the light ray transform $L$ on globally hyperbolic Lorentzian manifold $(\mcm, g)$ of dimension $n+1, n\geq 2$ which is null-geodesic complete and without conjugate points. Then $L: H^s_{\comp}(\mcm) \rightarrow H_{\loc}^{s+ s_0/2}(\mcc)$ is continuous with $s_0$ such that $\max(-n/2 + 1/2, -1) \leq -s_0$, $n\geq 2.$ 
\end{theorem}
 
Next, we use the microlocal structure of the normal operator to answer some questions regarding the determination of singularities. Consider the cotangent bundle $T^*\mcm$. Our convention for the signature of the metric $g$ is $(-, +, \cdots, +)$. A covector $\zeta\in T_z^*\mcm$ is called space-like if $g(\zeta, \zeta) > 0$, time-like if $g(\zeta, \zeta) < 0$ and light-like if $g(\zeta, \zeta) = 0$. The set of space-like, time-like and light-like vectors are denoted by $\Gamma^{sp}, \Gamma^{tm}$ and $\Gamma^{lt}$ respectively.  In general relativity, space-like singularities corresponds to particles moving slower than the speed of light, and light-like singularities corresponds to objects moving at the speed of light such as photons and gravitational waves. The following statements can be seen  from Theorem \ref{thm-main}, see also  \cite{LOSU2, StUh, Gui}. Suppose $f\in \mce'(\mcm)$:
\begin{enumerate}
\item If $\WF(f)\subseteq \Gamma^{sp}$, then $q\in \WF(Nf)$   if and only if $q\in \WF(f)$. 
\item If $\WF(f)\subseteq \Gamma^{tm}$, then $\WF(Nf) = \emptyset$. 
\end{enumerate}
It is interesting to know what happens to the light-like singularities which is currently unclear. We prove
\begin{theorem}\label{thm-main1}
For $(\mcm, g)$ as in Theorem \ref{thm-main} with $n=2, 3$, there exists $f\in \mce'(\mcm)$ such that $\WF(f) \subseteq \Gamma^{lt}$ and $\WF(f) \neq \emptyset$ but $\WF(Nf) = \emptyset$.
\end{theorem}

For such $f$, we have $Nf\in C^\infty(\mcm)$ and $Lf \in C^\infty(\mcc)$. The result implies that one may not be able to determine light-like singularities of $f$ using singularities of $Nf$. We remark that related examples of non-compactly supported distributions are known for $\mbr^{2+1}$ with the Minkowski metric, see \cite[Section 2]{GrUh}. Also, at least for the Minkowski spacetime, the transform is known to be injective on $C_0^\infty$. However, the proof is based on analyticity of the Fourier transform of $f$. Our result indicates that the determination is not stable when $f$ has light-like singularities.  

Under stronger conditions for example if the singularities of $f$ are of conormal type with principal symbols of a fixed sign, we prove in Theorem \ref{thm-main2} that the wave front set of $f$ can be determined from $Nf$. Another related example was found in \cite{VaWa} for small perturbations of Minkowski metric in the study of the inverse Sachs-Wolfe problem. In \cite{VaWa}, the author and A. Vasy analyzed the singularities in $Nf$ where $f$ is a solution to the Cauchy problem of wave equations. We proved a stronger statement that  $f$ can be stably determined by $Lf$. There has been active researches on the light ray transform recently,  see for example \cite{FIO, Ilm, LOSU2, Ste}. 

In this work, we mainly analyze the light ray transform acting on compactly supported distributions and  on Lorentzian manifolds without  conjugate points. 
For some applications in general relativity, conjugate points cannot be avoided, see for example the singularity theorem of Hawking and Penrose, see \cite{HaEl}.  In Section \ref{sec-conjugate}, we analyze the structure of the normal operator for standard static space-times with time-like conjugate points and show that the Schwartz kernel is the sum of a paired Lagrangian distribution as in Theorem \ref{thm-main} and a Lagrangian distribution associated with the conjugate points. One can derive the analogue of Theorem \ref{thm-main1}. Here, we benefited from the work of Stefanov and Uhlmann \cite{StUh} for the geodesic ray transform with caustics of fold type in the Riemannian setting.

The paper is organized as follows. In Section \ref{sec-pair} and \ref{sec-minres}, we prove the main results for the Minkowski spacetimes. These discussions are meant to serve as an introduction to the method and some microlocal machinery. Then we prove the theorem for globally hyperbolic manifolds in Section \ref{sec-sym} and Section \ref{sec-mainres}. We discuss the determination of light-like singularities in Section \ref{sec-det}. In Section \ref{sec-conjugate}, we address the case of standard static spacetime with time-like conjugate points.

\section{Preliminaries}\label{sec-pair}
It is beneficial to start with the light ray transform on $n+1$-dimensional Minkowski space  $(\mbr^{n+1}, g), n \geq 2$ where $g = -dt^2 + dx_1^2 + \cdots + dx_n^2$. Hereafter, we use $(t, x_1, \cdots, x_n)$ for the coordinates on $\mbr^{n+1}$. In fact, some of the calculations will be used for the general case later. We parametrize the future pointing light-like geodesics as follows: for $\theta\in \mbs^{n - 1}, z\in \mbr^n$, the light-like geodesics from $(0, z)$ in the direction $(1, \theta)$ is given by $l_{z, \theta}(s) = (s, z + s\theta), s\in \mbr.$ 
The set of light rays $\mcc$ is identified with $\mbr^n\times \mbs^{n-1}$ with the standard product measure.  
Using this parametrization,   the light ray transform is 
\beqq\label{eq-minlray}
Lf(z, \theta) = \int_\mbr f(s, z + s\theta) ds, \quad f\in C_0^\infty.
\eeqq 
Let $L^t$ be the transpose of $L$. Consider the normal operator $N = L^t L$. It is computed in \cite[Theorem 2.1]{LOSU2} that 
\beqq\label{eq-minnormal}
N f(t, x)  = \int_{\mbr^{n+1}} K_N(t, x, t', x') f(t', x') dt'dx'
\eeqq
where the Schwartz kernel  
\beqq\label{eq-norker}
K_N(t, x, t', x') = \frac{\delta(t - t' - |x - x'|) + \delta(t - t' + |x - x'|)}{|x - x'|^{n - 1}}
\eeqq 
In particular, $N$ can be written as an Fourier multiplier 
\beqq\label{eq-normalfourier}
\begin{gathered}
N f(t, x) =   \int_{\mbr^{n+1}} e^{i(t \tau +x\cdot \xi)} k(\tau, \xi) \hat f(\tau, \xi) d\tau d\xi
\end{gathered}
\eeqq
where 
\beqq\label{eq-ksym}
k(\tau, \xi) = C_n \frac{(|\xi|^2 - \tau^2)_+^{\frac{n - 3}{2}}}{|\xi|^{n-2}}, \quad C_n = 2\pi |\mbs^{n-2}|.
\eeqq
Here, for $s\in \mbr$,  $s_+^a, \re a> -1$ denotes the  distribution defined by $s_+^a = s^a$ if $s > 0$ and $s_+^a = 0$ if $s \leq 0$, see \cite[Section 3.2]{Ho1}. We can write the normal operator as
\beqq\label{eq-Nf}
N  f(t, x) =   \int_{\mbr^{n+1}} \int_{\mbr^{n+1}} e^{i(t - t', x - x')\cdot (\tau, \xi)} k(\tau, \xi)  f(t', x') d\tau d\xi dt'dx'
\eeqq
   
On the dual space $\mbr^{n+1}_{(\tau, \xi)}$, we let $\Gamma^{tm}_\pm = \{(\tau, \xi)\in \mbr^{n+1}: \tau^2 > |\xi|^2, \pm \tau > 0\}$ be the set of future/past pointing time-like covectors, and $\Gamma^{tm} = \Gamma^{tm}_+\cup \Gamma^{tm}_-$. Let $\Gamma^{sp} =  \{(\tau, \xi)\in \mbr^{n+1}: \tau^2 < |\xi|^2\}$ be the set of space-like  covectors. Finally, let $\Gamma^{lt}_\pm = \{(\tau, \xi)\in \mbr^{n+1}: \tau^2 = |\xi|^2, \pm \xi_0 > 0\}$ be the set of future/past pointing light-like  covectors. We also let $\Gamma^{lt} = \Gamma_+^{lt}\cup \Gamma_-^{lt}$. 
We see that in \eqref{eq-ksym}, the symbol $k(\tau, \xi)$ is supported in $\Gamma^{sp}$,  is homogeneous of degree $-1$ in $(\tau, \xi)$ and smooth away from $\Gamma^{lt}$. Moreover,  $k(\tau, \xi) \sim \text{dist}((\tau, \xi), \Gamma^{lt})^{(n - 3)/2},$ for $(\tau, \xi)$ space-like near $\Gamma^{lt}.$ Therefore, $k(\tau, \xi)$ looks like a symbol for a pseudo-differential operator of order $-1$ with a conormal singularity at $\Gamma^{lt}$. This is a typical example of operator whose Schwartz kernel is a paired Lagrangian distribution introduced in \cite{GuUh}. We recall some results on paired Lagrangian distributions that we use later. Our main reference is  Section 5 of \cite{DUV}.  

Let $\mcx$ be a $C^\infty$ manifold of dimension $n$ and $w_\mcx$ be the simplectic form on $T^*\mcx$. 
Let $\La_0, \La_1$ be conic Lagrangian submanifolds of $T^*(\mcx\times \mcx)\backslash 0$ with symplectic form $\pi_1^*w_\mcx + \pi_2^*w_\mcx$. Here, $\pi_1, \pi_2: \mcx\times \mcx\rightarrow \mcx$ denotes the projections to the first, second copy of $\mcx$. 
Suppose that $\La_1$ intersects $\La_0$ cleanly at a codimension $k$, $1\leq k\leq 2n-1$ submanifold $\Sigma = \La_0\cap \La_1$, namely
\beq
T_p(\La_0\cap \La_1) = T_p(\La_0)\cap T_p(\La), \quad \forall p\in \Sigma. 
\eeq
From \cite[Proposition 2.1]{GuUh}, we know that all such intersecting pairs $(\La_0, \La_1)$ are locally symplectic diffeomorphic to each other. To define paired Lagrangian distributions, we first consider the following model problem. 

Let $\tilde \mcx = \mbr^n = \mbr^k\times \mbr^{n-k}, 1\leq k\leq n-1$, and use coordinates $x = (x', x''), x'\in \mbr^k, x''\in \mbr^{n-k}$. Let $\tilde \La_0 = \{(x, \xi, x, -\xi)\in T^*(\tilde \mcx\times \tilde \mcx)\backslash 0 : \xi\neq 0\}$ be the punctured conormal bundle of $\diag$ in $T^*(\tilde \mcx \times \tilde \mcx)$, and 
\beq
\tilde \La_1 = \{(x, \xi, y, \eta)\in T^*(\tilde \mcx\times \tilde \mcx)\backslash 0: x'' = y'', \xi' = \eta' = 0, \xi'' = \eta'' \neq 0\} 
\eeq
which is the punctured conormal bundle to $\{(x, y)\in \tilde \mcx\times \tilde \mcx: x'' = y''\}$.  The two Lagrangians intersect cleanly at $\tilde \Sigma = \{(x, \xi, y, \eta)\in T^*(\tilde \mcx\times \tilde \mcx)\backslash 0 : x'' = y'', \xi'' = \eta'', x' = y', \xi' = \eta' = 0\}$ which is of codimension $k.$ For this model pair, the paired Lagrangian distribution $I^{p, l}(\mbr^n\times \mbr^n; \tilde \La_0, \tilde \La_1)$ consists of oscillatory integrals 
\beqq\label{eq-upair}
u(x, y) = \int e^{i[(x' - y' - s)\cdot \eta' + (x'' - y'')\cdot \eta'' + s\cdot \sigma]}a(s, x, y, \eta, \sigma) d\eta d\sigma ds
\eeqq
where $a$ is a product type symbol  which is a $C^\infty$ function and satisfies 
\beqq\label{eq-symord}
|\p_\eta^\alpha \p_\sigma^\beta \p_s^\theta \p^\gamma_x \p^\delta_y a(s, x, y, \eta, \sigma)| \leq C (1 + |\eta|)^{p + k/2 -|\alpha|}(1 + |\sigma|)^{l   - k/2 - |\beta|}
\eeqq
for multi-indices $\alpha, \beta, \theta, \gamma, \delta$ over each compact set $\mck$ of $\mbr^n\times \mbr^n \times \mbr^k.$ The constant $C$ depends on the indices and $\mck.$ The set of product type symbols is denoted by $S^{p, l}(\mbr^n\times \mbr^n; \mbr^n; \mbr^k)$. 
In \eqref{eq-symord}, the order of the symbol is chosen such that away from $\tilde \Sigma = \tilde \La_0\cap \tilde \La_1$, $u \in I^{p+l}(\mbr^n\times \mbr^n; \tilde \La_0)$ and  $u \in I^{p}(\mbr^n\times \mbr^n; \tilde \La_1)$ using H\"ormander's notion of Lagrangian distributions, see \cite[Section 25.1]{Ho4}. More precisely, for $u$ defined in \eqref{eq-upair},  the principal symbols are
\beq
\begin{gathered}
 a(0, y, y, \eta, \eta') |d\eta dy|^\ha \text{ on $\tilde \La_0\backslash \tilde \Sigma$}\\
( \int a(x' - y', x, y, 0, \eta'', \sigma) e^{i(x' - y')\sigma}d\sigma) |dy' d\eta''|^\ha \text{ on $\tilde \La_1\backslash \tilde \Sigma$}
\end{gathered}
\eeq
respectively. Here, the Maslov factors are not shown. We see that if $u \in I^p(\mbr^n\times \mbr^n; \tilde \La_1)$, the order $\mu$ of the symbol should be  
\beq
\mu - (2n)/4 + (n-k)/2 = p \Longrightarrow \mu = p + k/2
\eeq 
Then one can determine $l$ from that the order of the pseudo-differential operator symbol is $p+l.$  
 We also use the notation $I^{p, l}(\mbr^n\times \mbr^n; \tilde \La_0, \tilde \La_1)$ to denote the space of operators $A: \mce'(\mbr^n; \Omega^\ha_{\mbr^n}) \rightarrow \mcd'(\mbr^n; \Omega^\ha_{\mbr^n})$ where $\Omega^\ha_{\mbr^n}$ denotes the line bundle of half-densities on $\mbr^n$,  whose Schwartz kernel $K_A$ is a paired Lagrangian distribution with values in $\Omega_{\mbr^n\times \mbr^n}^\ha$. 
 
 There is an equivalent description of the paired Lagrangian distribution \eqref{eq-upair} introduced in \cite[Section 5]{DUV} that is convenient sometimes. Modulo $C_0^\infty(\mbr^n\times \mbr^n)$, \eqref{eq-upair} can be written as 
 \beqq\label{eq-upair1}
u(x, y) = \int e^{i[(x' - y')\cdot \eta' + (x'' - y'')\cdot \eta'']}b(x, y, \eta) d\eta 
\eeqq
where $b$ satisfies the following estimates. First, in the region $|\eta'|\leq C |\eta''|, |\eta''|\geq 1$, $b$ satisfies
\beq
|(Qb)(x, y, \eta)|\leq C\langle \eta''\rangle^{p+k/2} \langle \eta'\rangle^{l - k/2}
\eeq
for all $Q$ which is a finite product of differential operators of the form $D_{\eta'}, \eta'_j D_{\eta'_m}, \eta''_j D_{\eta_m''}$. Second, in the region $|\eta''|\leq C |\eta'|, |\eta'|\geq 1$, $b$ satisfies the standard regularity estimate 
\beq
|(Qb)(x, y, \eta)|\leq C\langle \eta'\rangle^{p+ l}
\eeq
for all $Q$ which is a finite product of differential operators of the form $\eta_j'D_{\eta'_m}, \eta'_j D_{\eta''_m}.$ We refer the readers to \cite[Section 5]{DUV} for the argument of the equivalence of \eqref{eq-upair} and \eqref{eq-upair1}, and a description of the principal symbols \cite[Lemma 5.3]{DUV}. 

For general cleanly intersecting pairs $(\La_0, \La_1),$ let $\chi: T^*(\mcx\times \mcx)\backslash 0\rightarrow T^*(\tilde \mcx\times \tilde \mcx)\backslash 0$ be a canonical transformation such that $\chi(\La_0) \subseteq \tilde \La_0, \chi(\La_1)\subseteq \tilde \La_1$. Then the set of paired Lagrangian distributions $I^{p, l}(\mcx \times \mcx; \La_0, \La_1)$ are defined invariantly by conjugating elements of $I^{p, l}(\tilde \mbr^n\times \tilde \mbr^n; \tilde \La_0, \tilde \La_1)$ by Fourier integral operators with canonical relation $\chi$, see \cite{GuUh} for more details. It is known that for any $u\in I^{p, l}(\mcx\times \mcx; \La_0, \La_1)$, $\WF(u)\subseteq \La_0 \cup \La_1$. One can deduce that on $\La_0\backslash \Sigma$, $u$ is a pseudo-differential operator of order $p+l$. On $\La_1\backslash \Sigma$, $u\in I^p(\mcx \times \mcx; \La_1)$ is a Fourier integral operator of order $p$. From \cite[Proposition 6.2]{GuUh}, we know that 
\beq
\bigcap_{l \in \mbr}I^{p, l}(\mcx\times \mcx; \La_0, \La_1) = I^p(\mcx\times \mcx; \La_1), \quad \bigcap_{p\in \mbr}I^{p, l}(\mcx\times \mcx; \La_0, \La_1) = C^\infty(\mcx).
\eeq

To define symbols of distribution $u\in I^{p, l}(\mcx\times \mcx; \La_0, \La_1)$, we note that there are two well-defined principal symbols (valued in half-densities tensored with Maslov bundles) on $\La_0\backslash \Sigma$, denoted by $\sigma_0(u)$, and on $\La_1\backslash \Sigma$, denoted by $\sigma_1(u)$. These symbols are singular at $\Sigma$ 
and the space of such such symbols is denoted by $S^{p, l}(\mcx \times \mcx; \La_0, \La_1)$   together with the line bundle of half densities tensored with the Maslov bundle, see \cite[Section 5]{GrUh0}. We remark that $\sigma_0(u)$ determines the leading order terms of $\sigma_1(u)$ at $\Sigma$, however, it does not determine the full $\sigma_1(u)$ on $\La_1\backslash \Sigma.$ 
 We will also need the Sobolev estimates for paired Lagrangian distributions for the flow out model. A submanifold $\Gamma \subseteq T^*\mcx$ is involutive if $\Gamma = \{(x, \xi): p_i(x, \xi) = 0, i = 1, 2, \cdots, k\}$ satisfies (i) $p_i$ are defining functions of $\Gamma$  that $dp\neq 0$ on $p = 0$,  and (ii) $p_i$ are in involution so the Poisson brackets $\{p_i, p_j\} = 0$ at $\Gamma$. Let $H_{p_i}$ be the Hamilton vector fields of $p_i$. The flow out of $\Gamma$ 
\beqq\label{eq-flowout}
\begin{gathered}
\La_\Gamma = \{(x, \xi, y, \eta)\in T^*\mcx \times T^*\mcx: (x, \xi)\in \Gamma,  
(y, \eta) = \exp(\sum_{j = 1}^k t_j H_{p_j})(x, \xi), t_j \in \mbr \}
\end{gathered}
\eeqq
is a Lagranian submanifold of $T^*(\mcx \times \mcx)$ and is a canonical relation if $\Gamma$ is conic.
 \begin{theorem}[Proposition 5.6 of \cite{DUV}]\label{thm-sobolev}
 Let $A\in I^{p, l}(\mbr^n\times \mbr^n; \La_0, \La_1)$ with $\La_1$ the flow out $\La_\Sigma$. Then $A: H^s_{\comp}(\mbr^n)\rightarrow H^{s+ s_0}_{\loc}(\mbr^n)$ is continuous for any $s\in \mbr$ if
 \beq
 \max(p + k/2, p + l)\leq -s_0.
 \eeq
\end{theorem}

\section{Microlocal results for the Minkowski light ray transform}\label{sec-minres}
To see that the kernel \eqref{eq-norker} is a paired Lagrangian distribution, we proceed directly using the oscillatory integral representation \eqref{eq-Nf}. The Schwartz kernel of $N$ is 
\beq 
K_N(t, x, t', x') =   \int_{\mbr^{n+1}} e^{i(t - t')\tau+ i(x - x')\cdot \xi} k(\tau, \xi)  d\tau d\xi, \quad k(\tau, \xi) = C_n \frac{(|\xi|^2 - \tau^2)_+^{\frac{n - 3}{2}}}{|\xi|^{n-2}}.
\eeq 
The symbol $k$ has singularities at $\xi = 0$  but these can be removed by using compactly supported smooth cut-off functions, which changes $K_N$ by a $C^\infty$ term. The singularity at $\tau^2 - |\xi|^2 = 0$ can be removed similarly because the symbol is integrable there.  To simplify the notations, we will not spell out the cut-offs.   We let $s = \tau - |\xi|$ and write  
\beqq\label{eq-ks}
\begin{gathered}
K_N(t, x, t', x') =   \int_{\mbr^{n+1}} e^{i(t - t')(s + |\xi|) + i(x - x')\cdot \xi} k(s, \xi)  ds d\xi, \\
 k(s, \xi) = C_n \frac{s_-^{\frac{n - 3}{2}}(s +2|\xi|)_+^{\frac{n-3}{2}}}{|\xi|^{n-2}} +  C_n \frac{s_+^{\frac{n - 3}{2}}(s +2|\xi|)_-^{\frac{n-3}{2}}}{|\xi|^{n-2}} =  C_n \frac{s_-^{\frac{n - 3}{2}}(s +2|\xi|)_+^{\frac{n-3}{2}}}{|\xi|^{n-2}}.
\end{gathered}
\eeqq
To see that this can be transformed to the model integral \eqref{eq-upair1}, consider the simplectic change of variables  on $T^*\mbr^{n+1}$
\beq
\tilde x = x - (t - t')\xi/|\xi|, \quad \tilde t = t - t', \quad s = s, \quad \xi = \xi.
\eeq
We can choose an Fourier integral operator with symbol of order $0$ which quantizes the symplectic change of variable to transform $K_N$ to
\beqq\label{eq-kker}
\begin{gathered}
K_N(\tilde t, \tilde x, t', x')   =  \int_{\mbr^{n+1}} e^{i\tilde t s + i \tilde x\cdot \xi} k(s, \xi) ds d\xi 
 \end{gathered}
\eeqq 
modulo a smooth term. 
The symbol $k(s, \xi)$ satisfies the product type estimate with $p = -n/2, l = n/2-1$. In fact, for $|\xi| \leq C|s|, |s|\geq 1$, we have 
\beq
|k(s, \xi)| \leq C \frac{|s|^{\frac{n - 3}{2}}|s|^{\frac{n-3}{2}}}{|\xi|^{n-2}} \leq C|s|^{-1}
\eeq
One can verify the same estimate for $Qk$ where $Q$ is the finite product of differential operators of the form 
$sD_{s}, s D_{\xi_m}, m = 1, 2, \cdots, n$. For $|s|\leq C|\xi|, |\xi|\geq 1$, we have 
\beq
|k(s, \xi)| \leq C\frac{|s|^{\frac{n - 3}{2}} |\xi|^{\frac{n-3}{2}}}{|\xi|^{n-2}} \leq C  |\xi|^{-n/2 + 1/2} |s|^{n/2-3/2}
\eeq
 and one can verify the estimate for $Qk$ where $Q$ is the finite product of differential operators of the form $D_s, s D_{s}, \xi_j D_{\xi_m}$. So $K_N$ is a paired Lagrangian distribution. The two associated Lagrangians are  
 \beqq\label{eq-lag1}
 \La_0 = \{(t, x, \tau, \xi; t', x', \tau', \xi')\in T^*\mbr^{n+1}\backslash 0 \times T^*\mbr^{n+1}\backslash 0: t'  = t, x' = x, \tau' = -\tau, \xi' = -\xi\}
 \eeqq
which is the punctured conormal bundle of the diagonal in $\mbr^{n+1}\times \mbr^{n+1}$ and 
 \beqq\label{eq-lag2}
 \begin{gathered}
 \La_1 = \{(t, x, \tau, \xi; t', x', \tau', \xi')\in T^*\mbr^{n+1}\backslash 0 \times T^*\mbr^{n+1}\backslash 0: x = x' + (t - t')\xi/|\xi|,\\
 \tau = \pm|\xi|,  \tau' = -\tau,  \xi' = -\xi \}.
  \end{gathered}
 \eeqq
Let $f(\tau, \xi) = -\tau^2 - |\xi|^2$ and $\Sigma = \{(t, x, \tau, \xi; t, x, -\tau, -\xi)\in T^*\mbr^{n+1}\backslash 0 \times T^*\mbr^{n+1}\backslash 0: f(\tau, \xi) = 0\}$.  Then $\La_1$ is the flow out of $\Sigma$ under the Hamilton vector field $H_f$ of $f(\tau, \xi) = \ha(\tau^2 - |\xi|^2).$ Indeed, we have 
\beq
H_f = \tau \frac{\p}{\p t} + \sum_{i = 1}^3 \xi_i \frac{\p}{\p x^i}
\eeq
Let $\gamma(s) = (t(s), x(s), \tau(s), \xi(s))$ be a null bi-characteristic which satisfy 
\beq
\begin{gathered}
\dot t(s) = \tau, \quad \dot x_i(s) = \xi_i, \quad \dot \tau(s) = 0, \quad \dot \xi_i(s) = 0\quad s\geq 0\\
 t(0) = t', \quad x(0) = x', \quad \tau(0) = \tau', \quad \xi(0) = \xi'
\end{gathered}
\eeq
with $f(\tau', \xi') = 0.$ We solve that 
\beq
t(s) = t' +  s \tau', \quad x(s) = x' + s \xi', \quad \tau(s) = \tau', \quad \xi(s) = \xi'.  
\eeq
which up to a re-parametrization gives \eqref{eq-lag2}. Thus $\La_1$ intersects $\La_0$ cleanly at $\Sigma.$

From \eqref{eq-Nf}, we know that $N$ on $\La_0\backslash \La_1$ is a pseudo-differential operator of order $(n-3) - (n-2) = -1$.  The principal symbol of $N$ on $\La_0\backslash \La_1$ can be read from \eqref{eq-Nf} and it is non-vanishing.  
Next, we find the oscillatory integral representation of $N$ on $\La_1\backslash \La_0$ and determine the principal symbol.  
We first use the kernel representation \eqref{eq-ks} and write it as
\beqq\label{eq-kker0}
\begin{gathered}
K_N(t, x, t', x')    
   = \int_{\mbr^{n}} e^{i(t - t') |\xi| + i(x - x')\cdot \xi} A(t - t', \xi)  d\xi 
 \end{gathered}
\eeqq  
where 
\beqq\label{eq-syma}
\begin{gathered}
A(\sigma, \xi) = \int_{-2|\xi|}^{0}e^{i\sigma s} k(s, \xi) ds  = \int_{-1}^0 e^{i\sigma 2|\xi| s} k(2 |\xi| s, \xi) 2 |\xi| ds  \\
=  2^{n-2} C_n \int_{0}^1 e^{-i\sigma 2|\xi| s}  s^{\frac{n - 3}{2}} (1 - s)^{\frac{n-3}{2}} ds 
\end{gathered}
\eeqq
For $n= 3$, 
\beq
\begin{gathered}
A(\sigma, \xi)  
=  2 C_3 \int_{0}^1 e^{-i\sigma 2|\xi| s}   ds   = C_3 \frac{1}{i \sigma |\xi|} (1-  e^{-2i\sigma |\xi|}) 
\end{gathered}
\eeq
Then $K_N$ becomes 
\beqq\label{eq-kn3} 
\begin{gathered}
K_N(t, x, t', x')    =   C_3 \int_{\mbr^{n}} e^{i(t - t') |\xi| + i(x - x')\cdot \xi}  \frac{1}{i (t-t') |\xi|} (1 - e^{-2i(t - t')|\xi|})  d\xi \\
 = C_3 \int_{\mbr^{n}} e^{i(t - t') |\xi| + i(x - x')\cdot \xi}  \frac{1}{i (t-t') |\xi|}  d\xi  + C_3 \int_{\mbr^{n}} e^{-i(t - t') |\xi| + i(x - x')\cdot \xi}  \frac{1}{i (t-t') |\xi|}   d\xi 
 \end{gathered}
\eeqq 
We see that this is an Fourier integral operator associated with $\La_1$ with symbol of order $-1$ when $t\neq t'.$  

For $n\geq 5$ odd, we use \eqref{eq-syma} and apply integration by parts  to get
\beq
\begin{gathered}
A(\sigma, \xi)  =  
 C_n \frac{2^{n-3}}{i\sigma|\xi|}  
\int_{0}^1 e^{-i\sigma 2|\xi| s}  (s^{\frac{n - 3}{2}}(1 - s)^{\frac{n-3}{2}})' ds \\
 = C_n \frac{2^{n-3}}{i\sigma|\xi|}  
\int_{0}^1 e^{-i\sigma 2|\xi| s} \frac{n-3}{2} (s^{\frac{n - 3}{2}-1}(1 - s)^{\frac{n-3}{2}} - s^{\frac{n - 3}{2}}(1 - s)^{\frac{n-3}{2}-1}) ds
\end{gathered}
\eeq
Repeating the integration by part $\frac{n-3}{2}$ times, we get 
\beq
\begin{gathered}
A(\sigma, \xi)   = C_n \frac{1}{(2i \sigma |\xi|)^{\frac{n-3}{2}}} (\frac{n-3}{2})! \int_0^1 e^{-i\sigma 2|\xi|s}  (1 - s)^{(n-3)/2}ds \\
 + C_n \frac{1}{(2i \sigma |\xi|)^{\frac{n-3}{2}}} (\frac{n-3}{2})! (-1)^{(n-3)/2}\int_0^1 e^{-i\sigma 2|\xi|s} s^{(n-3)/2}ds \\
 + \sum_{k, j \geq 1, k + j = (n-3)/2} c_{k, j} \int_0^1 e^{-i\sigma 2|\xi|s}  s^{(n-3)/2-k}(1 - s)^{(n-3)/2-j}ds
\end{gathered}
\eeq
where $c_{k, j}$ are constants. So far, the boundary terms from integration by parts vanish. We continue with integration by parts to get 
\beq
\begin{gathered}
A(\sigma, \xi)   = C_n \frac{1}{(2i \sigma |\xi|)^{\frac{n-3}{2}} (2i\sigma |\xi|)} (\frac{n-3}{2})! + \sum_{k = 1}^M a_k(\sigma ) |\xi|^{-\frac{n-3}{2} -1 - k} \\
 + C_n \frac{1}{(2i \sigma |\xi|)^{\frac{n-3}{2}}(-2i\sigma |\xi|)} (\frac{n-3}{2})! (-1)^{(n-3)/2} e^{-i\sigma 2|\xi|} + \sum_{k = 1}^M b_k(\sigma ) e^{-i\sigma 2|\xi|}  |\xi|^{-\frac{n-3}{2} -1 - k}  \end{gathered}
\eeq
where $a_k, b_k$ are smooth in $\sigma$ for $\sigma \neq 0$, and $M$ is some constant depending on $n. $ We remark that for $n$ even, integration by parts will eventually lead to singular integrals. This is why we cannot deal with $n$ even at this point. We let 
\beq
\begin{gathered}
a_0(\sigma) = C_n \frac{1}{(2i \sigma)^{\frac{n-3}{2} + 1}}  (\frac{n-3}{2})!, \quad 
b_0(\sigma) = C_n \frac{1}{(2i \sigma)^{\frac{n-3}{2}+1}}  (\frac{n-3}{2})! (-1)^{(n-3)/2+ 1}  
\end{gathered}
\eeq
 We see that 
\beq
\begin{gathered}
A(\sigma, \xi)   =  \sum_{k = 0}^M a_k(\sigma ) |\xi|^{-\frac{n-3}{2} -1 - k}  + \sum_{k = 0}^M b_k(\sigma ) e^{-i\sigma 2|\xi|}  |\xi|^{-\frac{n-3}{2} -1 - k}  
\end{gathered}
\eeq
Finally, using \eqref{eq-kker0}, we get 
\beqq\label{eq-ksym1} 
\begin{gathered}
K_N(t, x, t', x') 
 =  \int_{\mbr^{n}} e^{i(t - t') |\xi| + i(x - x')\cdot \xi}  \sum_{k = 0}^M a_k(\sigma ) |\xi|^{-\frac{n-3}{2} -1 - k} d\xi  \\
 + \int_{\mbr^{n}} e^{-i(t - t') |\xi| + i(x - x')\cdot \xi} \sum_{k = 0}^M b_k(\sigma )   |\xi|^{-\frac{n-3}{2} -1 - k}  d\xi 
 \end{gathered}
\eeqq 
which shows that $K_N$ is an Lagrangian distribution associated with $\La_1$ with a symbol of order $-(n-3)/2 - 1$. Thus, $N$ on $\La_1\backslash \La_0$ is an FIO of order $p =  -n/2+1/2 - (n+1)/2 + n/2 = -n/2$. Because $p+l = -1$, we get $l = n/2 - 1.$ The principal symbol is seen non-vanishing from \eqref{eq-ksym1}.

To deal with $n$ even, we will use another representation of the kernel $K_N$ which leads to a very simple proof that $N$ is an elliptic FIO on $\La_1\backslash \La_0.$ For this, we go back to the kernel \eqref{eq-norker} and write it as 
\beqq\label{eq-norker}
K_N(t, x, t', x') = \int_\mbr e^{i((t - t') - |x - x'|)\tau} (t - t')_+^{-n + 1}d\tau  +\int_\mbr  e^{i(t - t' +|x - x'|)\tau} (t - t')_-^{-n + 1} d\tau  
\eeqq 
For $t \neq t'$, $K_N$ is a Lagrangian distribution conormal to the light cone  minus the vertex $\{(t, x, t', x')\in \mcm\times \mcm: t - t' = |x - x'|, t\neq t'\}$. Also, the principal symbol  is seen non-vanishing.   

To summarize, we proved the Minkowski version of Theorem \ref{thm-main}. 
\begin{theorem}\label{thm-mainmin0}
For the Minkowski light ray transform $L$ defined in \eqref{eq-minlray} and $n\geq 2$, the Schwartz kernel of the normal operator $K_N \in I^{-n/2, n/2- 1}(\mbr^{n+1}\times \mbr^{n+1}; \La_0, \La_1)$, in which $\La_0, \La_1$ are two cleanly intersection Lagrangians defined in \eqref{eq-lag1}, \eqref{eq-lag2}. Let $\Sigma = \La_0\cap \La_1.$ The principal symbols  of $K_N$ on $\La_0\backslash \Sigma, \La_1\backslash \Sigma$ are non-vanishing. 
\end{theorem}

Using the Sobolev estimate Theorem \ref{thm-sobolev} and the above proposition, we immediately obtain
\begin{cor}
The Minkowski light ray transform $L: H^s_{\comp}(\mbr^{n+1}) \rightarrow H_{\loc}^{s+ s_0/2}(\mbr^n\times \mbs^{n-1})$ is continuous  with $s_0$ such that $\max(-n/2 +1/2, -1) \leq -s_0$ 
for $n\geq 2$.
\end{cor}

Finally, we mention that it is possible to find a parametrix of $N$ in the same spirit as in \cite{GrUh0} and \cite{Wan}. The parametrix is a paired Lagrangian distribution and there is an error term belonging to $I^{-1/2}(\mbr^n; \La_1)$. We will not prove it here. Instead, we discuss the case for $n = 2, 3$ for which the error term can be removed.  
 To state our result, we introduce the following substitute of the identity operator. Let $\chi_{sp}$ be the characterisitic function of $\Gamma^{sp}$ in $\mbr_{(\tau, \xi)}^{n+1}$, we define 
 \beqq\label{eq-id}
H f(t,  x) = \int_{\mbr^{n+1}} \int_{\mbr^{n+1}} e^{i(t - t')\tau + i (x - x')\xi} \chi_{sp}(\tau, \xi)  f(t', x') d\tau d\xi dt' dx'
 \eeqq
The symbol has Heaviside type singularities at $\tau = \pm |\xi|$. By the same arguments in Theorem \ref{thm-mainmin0}, the Schwartz kernel of $H$ is a paired Lagrangian distribution associated with the pair $(\La_0, \La_1).$ 
In fact on $\La_0\backslash \La_1,$ $H$ is a pseudo differential operator  of order $0$ so $p + l = 0$. On $\La_1\backslash \La_0$, $H$ is an FIO with a symbol of order $0$. Actually, as in \eqref{eq-kker0} and \eqref{eq-syma}, the Schwartz kernel of $H$ is
\beqq\label{eq-hker} 
\begin{gathered}
K_H(t, x, t', x')   
   = \int_{\mbr^{n}} e^{i(t - t') |\xi| + i(x - x')\cdot \xi} A(t - t', \xi)  d\xi 
 \end{gathered}
\eeqq  
where 
\beqq\label{eq-syma1}
\begin{gathered}
A(\sigma, \xi) = \int_{-2|\xi|}^{0}e^{i\sigma s} \chi_{sp}(s, \xi) ds  = \int_{-1}^0 e^{i\sigma 2|\xi| s} \chi_{sp}(2 |\xi| s, \xi) 2 |\xi| ds  \\
=   2  C_n \int_{-1}^0 e^{i\sigma 2|\xi| s}|\xi|  ds  
=    C_n \frac{1}{i\sigma } (1 - e^{-i\sigma 2|\xi|})
\end{gathered}
\eeqq
for all $n\geq 2$. Thus 
\beq 
\begin{gathered}
K_H(t, x, t', x') 
 =  C_n \frac{1}{i(t - t')} [ \int_{\mbr^{n}} e^{i(t - t') |\xi| + i(x - x')\cdot \xi}  d\xi 
 + \int_{\mbr^{n}} e^{-i(t - t') |\xi| + i(x - x')\cdot \xi}   d\xi ] 
 \end{gathered}
\eeq 
So we have $p = 0 - (n+1)/2 + n/2 = -1/2$ and $l  = 1/2$. Thus $H \in I^{-1/2, 1/2}(\mbr^{n+1}\times \mbr^{n+1}; \La_0, \La_1)$. One can see that the  principal symbols on each Lagrangian  are non-vanishing.  
 We think of $H$ as the identity operator on the light cones up to the boundary of $\Gamma^{sp}$.

 \begin{prop}\label{prop-para1}
Consider the Minkowski light ray transform $L$ defined in \eqref{eq-minlray} for $n = 2, 3$.  
 There exists a paired Lagrangian distribution $Q \in I^{n/2-1, -n/2 +2}(\mbr^{n+1}\times \mbr^{n+1}; \La_0, \La_1)$ such that 
 \beq
Q\circ N = H
 \eeq 
 As a result, we obtain an operator $A = Q\circ L^t$ such that $A\circ L = H.$ 
 \end{prop}
 \bpf
 We let $Q$ be defined by 
 \beq
Q f(t, x) =   \int_{\mbr^{n+1}\times \mbr^{n+1}} e^{i(t - t')\tau +i(x - x')\cdot \xi}  (\tau^2 - |\xi|^2)^{-(n-3)/2}|\xi|^{n-2} f(t', x') d\tau d\xi dt'dx'
\eeq
One can check that $Q\circ N = \id$ using the Fourier transform. The kernel of $Q$ is a paired Lagrangian distribution for $n = 2, 3$. Actually, the symbol is singular at $\tau^2 = |\xi^2|$ but the singularities can be removed by introducing a smooth cut-off function for $n = 2, 3$ which changes $Q$ by a smoothing operator. 
By the arguments in Theorem \ref{thm-mainmin0}, $Q \in I^{n/2 - 1}(\mbr^{n+1}\times \mbr^{n+1};  \La_1)$ and $Q\in I^{1}(\mbr^{n+1}\times \mbr^{n+1}; \La_0)$, thus $Q \in I^{n/2-1, -n/2+2}(\mbr^{n+1}\times \mbr^{n+1}; \La_0, \La_1)$ which completes the proof.
 \epf

\section{The Schwartz kernel}\label{sec-sym}
We consider  the light ray transform on globally hyperbolic Lorentzian manifold $(\mcm, g)$ described in \eqref{eq-metric}.  Let $\mcf$ be the set of all geodesics on $(\mcm, g)$, and  $\mcc$ be the set of  light-like geodesics, so  $\mcc\subseteq \mcf$. Provided there is no conjugate points on $(\mcm, g)$, $\mcf$ is a  $2n$ dimensional smooth manifold and $\mcc$ is a codimension one submanifold. We view the light ray transform as an operator $L: C_0^\infty(\mcm)\rightarrow C^\infty(\mcc)$. It is known that the Schwartz kernel $K_L$ is the delta distribution supported on the point-geodesic relation 
\beq
\mcz = \{(z, \gamma)\in \mcm \times \mcc: z  \in \gamma \}
\eeq
Therefore, $L$ is an Fourier integral operator and the kernel has conormal singularities to $\mcz$. The canonical relation can be described using Jacobi fields as in the Riemmanian setting, see \cite{GrUh0}. 
Consider $T^*\mcf$ which is equipped with a symplectic structure. Let $\gamma \in \mcf$ be a geodesic and $\mcj$ be the set of Jacobi fields tangent to $\dot\gamma$. So $\mcj$ is spanned by $\{\dot \gamma(s), s \dot \gamma(s)\}$. Let $\mcj^\perp$ be the orthogonal complements of $\mcj$ in $T^*_{\gamma(s)}\mcm$ with respect to the symplectic structure. Then we have  $T_\gamma^* \mcf \simeq \mcj^\perp,  T_\gamma \mcf \simeq \mcj^\perp$. 
Next, let $f(z, \zeta) = \ha g(z, \zeta), (z, \zeta)\in T^*\mcm$. We have $\mcc = \{\gamma \in \mcf : g(\dot\gamma, \dot\gamma)= 0\}$. Therefore, the conormal bundle $N^*_\gamma\mcc$ is spanned by $(-d_z f, d_\zeta f)$ and is identified as a  subspace of $\mcj^\perp$. Using these notions, we get that the canonical relation of $L$ or the twisted conormal bundle of $\mcz$ is 
\beq
C  = \{(y, \eta, \gamma_{z, \zeta}, \Gamma): y \in \gamma_{z, \zeta}(s), g(\eta, \dot \gamma_{z, \zeta}(s)) = 0, \Gamma = (-d_z f, d_\zeta f), s\in \mbr \}
\eeq 
Using H\"ormander's notion, the Schwartz kernel $K_L \in I^{-n/4}(\mcc\times \mcm; C')$.  

For globally hyperbolic manifolds $(\mcm, g)$, there is a natural choice of parametrization of $\mcc$. Instead of $(\mcm, g)$ given in \eqref{eq-metric}, we will be working with 
\beqq\label{eq-ghmetric}
\mcm = \mbr\times \mcn, \quad g = -dt^2 + h(t, x, dx)
\eeqq
which is conformal to \eqref{eq-metric} so that the $\beta >0 $ factor in \eqref{eq-metric} is absorbed to $h.$ The light ray transform differs by a smooth weight function (because $\beta$ is smooth) and will not change the analysis below. 
Let $\bar g = h(0, x, dx)$ be the Riemannian metric on $\mcn$ and 
\beqq
 S\mcn = \{(x, \theta)\in T\mcn: \bar g(\theta, \theta) = 1\}
\eeqq
be the unit sphere tangent bundle of $(\mcn, \bar g)$. For $(x, \theta)\in S\mcn$,  $v = (1, \theta)$ is a future pointing light-like vector at $(0, x)$. We denote by $\gamma_{x, \theta}$ the light-like geodesic from $(0, x)$ in direction $v$. Using the exponential map on $(\mcm, g)$, we also write $\gamma_{x, \theta}(s) = \exp_{(0, x)} s v,  s\in \mbr.$ 
We will identify $\mcc = S\mcn$ with the Liouville measure. We study the light ray transform \eqref{eq-lray0}  parametrized by 
\beqq\label{eq-lorlray}
L f(x, \theta) = \int_\mbr f(\gamma_{x, \theta}(s)) d s, \quad f\in C^\infty_0(\mcm), \quad (x, \theta)\in S\mcn.
\eeqq
 
Before we analyze the microlocal structure of the normal operator, it is useful to look at the  double fibration picture following the point of view of Guillemin \cite{Gui1} 
\beqq\label{eq-double}
\begin{tikzpicture}
  \matrix (m) [matrix of math nodes, row sep=2em,column sep=3em,minimum width=2em]
  {   &  C  & \\
     T^*\mcm \backslash 0 & &  T^*\mcc\backslash 0 \\};
 \path[-stealth]
    (m-1-2) edge node [left] {$\pi_\mcm$} (m-2-1) 
    (m-1-2) edge node [right] {$\pi_\mcc$} (m-2-3);
\end{tikzpicture}
\eeqq
If $\pi_\mcc$ is an injective immersion, then the composition $L^t L$ can be studied using Duistermaat and Guillemin's clean FIO calculus, see \cite{Gui1}. If $\xi$ in $\mcc$ is space-like, $\pi_\mcc$ is indeed injective, see \cite{LOSU2}. But it fails to be injective when $\xi$ is light-like. For $\mbr^{2+1}$,  Guillemin \cite{Gui} showed that $\pi_\mcc$ is a fold and $\pi_\mcm$ is a blow down, and the composition $L^t L$ fits in the framework developed in Greenleaf and Uhlmann \cite{GrUh0}. We will not take this approach here. Instead, we will  give a more direct proof by analyzing the Schwartz kernel. We start with static spacetimes, for which the calculation is similar to that for the Minkowski spacetime. Then we study globally hyperbolic spacetimes.

\subsection{Static spacetimes}
Let's recall the following definitions from \cite[Section 7.2]{SaWu}. Let $Z$ be a reference frame on $\mcm$, meaning that  $Z$ is a vector field each of whose integral curve is a future pointing time-like curve $\gamma: \mbr\rightarrow \mcm$ such that $g(\dot \gamma, \dot \gamma) = -1$ (called an observer). Then $Z$ is called stationary if there is a positive function $f$ on $\mcm$ such that $fZ$ is a Killing vector field. In addition, if $Z$ is irrotational, then $Z$ is called static. $\mcm$ is called stationary (static) if there is a stationary (static) reference frame $Z$ on $\mcm.$

Let $I$ be an open interval and $\mcn$ be a $n$ dimensional smooth manifold. The standard static spacetime is the manifold $I\times S$ with metric $g = -f(x)^2 dt^2 + h(x, dx)$ 
where $h$ is a Riemannian metric on $\mcn,$ see \cite[Page 360]{One}. In this case, the metric is static relative to $f\p_t.$ Locally, static spacetimes are isometric to the standard static spacetime. Examples of static spacetimes include Minkowski spacetime, and Roberston-Walker spacetimes up to a conformal transformation. It is known that static spacetime $(\mcm, g)$ is globally hyperbolic if and only if $(\mcn, h)$ is a  complete Riemannian manifold, see \cite{One}. 

For simplicity, we assume in this subsection that $(\mcm, g)$ is a standard static spacetime of the form 
\beqq\label{eq-staticmetric}
\mcm = \mbr\times \mcn, \quad g = -dt^2 + h(x, dx)
\eeqq
and there is no conjugate points on $(\mcm, g)$. Here, $h$ is a Riemannian metric on $\mcn$. In this case, light-like geodesics on $(\mcm, g)$ are lifts of geodesics on $(\mcn, h)$. More precisely, let $(x, \theta)\in S\mcn$ so that $h(\theta, \theta) = 1$. Then we have 
\beqq\label{eq-geo0}
\gamma_{x, \theta}(s) = \exp_{(0, x)}s(1, \theta) = (s, \exp^{h}_x(s\theta) )
\eeqq
where $\exp^{h}$ denotes the exponential map on $(\mcn, h)$. One can see this by examining the geodesic equation. It is convenient to look at the bicharacteristic flow in $T^*\mcm$. Let $\zeta = (\tau, \xi)$ be a covector in $T_{(t, x)}^*\mcm$. Consider the Hamiltonian function 
\beqq\label{eq-hamilton0}
f(t, x, \tau, \xi) = \ha g(\zeta, \zeta) = \ha (-\tau^2 + h(\xi, \xi))
\eeqq
The characteristic set is $\{f = 0\} = \Gamma^{lt}$. The Hamilton vector field is given by 
\beq
H_f =  \frac{\p f}{\p \tau} \frac{\p }{\p t} + \frac{\p f}{\p \xi} \frac{\p }{\p x}  - \frac{\p f}{\p t} \frac{\p }{\p \tau} - \frac{\p f}{\p x} \frac{\p }{\p \xi}
\eeq
Let $\alpha(s) = (t(s), x(s), \tau(s), \xi(s))$ be a bicharacteristic curve in $T^*\mcm$. We get using \eqref{eq-hamilton0} that 
\beqq\label{eq-bisys0}
\begin{gathered}
\frac{dt}{ds} =  -\tau, \quad  \frac{d\tau}{ds} = 0, \\
\frac{dx^j}{ds} = \ha \p_{\xi_j} h(\xi, \xi)\quad \frac{d\xi_j}{ds} =  -\ha\p_{x^j} h(\xi, \xi)
\end{gathered}
\eeqq
We see that the system is decoupled. Let's consider initial conditions
\beq
t(0) = t_0, \quad \tau(0) = \tau_0, \quad x(0) = x_0, \quad \xi(0) = \xi_0. 
\eeq
It is known that the projection of $\alpha(s)$ to $\mcm$ is a geodesic.  
We write $\zeta_0 = (\tau_0, \xi_0)$ and take $\tau_0 = 1, |\xi_0|_h = 1$ so that $\zeta_0$ is light-like, that is $\zeta_0\in \Gamma^{lt}$. We see that $\tau = 1$ is constant along the bicharacteristic.  More precisely,  $(t(s), x(s)) = \exp_{(t_0, x_0)}(s \zeta_0^\#)$. Using \eqref{eq-bisys0}, we see that $t(s) = s$ and $x(s) = \exp_{x_0}^h (s\xi_0^\#)$. So we arrive at the decomposition \eqref{eq-geo0}. 
The decomposition \eqref{eq-geo0} also implies that for static spacetime \eqref{eq-staticmetric}, $(\mcm, g)$ has no conjugate points if and only if  $(\mcn, h)$ has no conjugate points.  

Using $(x, \theta)\in S\mcn$ to parametrize the light rays, the light ray transform \eqref{eq-lorlray}  becomes 
\beqq\label{eq-staticlray}
L (f)(x, \theta) = \int_\mbr f(s, \exp^{h}_x(s \theta)) d s, \quad f\in C^\infty_0(\mcm).
\eeqq  
On $\mcm$, we will use the measure $dg = dt\wedge d\vol_h(x)$ which is $\sqrt{\det h}dt\wedge dx$ in local coordinates. 

\begin{prop}\label{prop-normal0}
For the light ray transform  \eqref{eq-staticlray} on a static spacetime \eqref{eq-staticmetric} of dimension $n+1, n\geq 2$ and without conjugate points, the Schwartz kernel $K_N$ of the normal operator $N = L^tL$  is 
\beqq\label{eq-normal0}
K_N(t, x, t', x') = \frac{\delta(t - t' - \dist^{h}(x, x')) + \delta(t - t' + \dist^{h}(x, x'))}{(\dist^{h}(x, x'))^{n-1}} \mcj(x, x')
\eeqq
for $(t, x), (t', x')\in \mcm.$ Here, $\dist^{h}: \mcn\times \mcn \rightarrow [0, \infty)$ is the distance function on $(\mcn, h)$ and $\mcj$ is a smooth non-vanishing function on $\mcn\times \mcn$ with $\mcj(x, x) = 1, x\in \mcn.$ 
\end{prop}
\bpf
We compute $L^t$. Using the Liouville measure $d\vol_h \wedge d\theta$ on $S\mcn$, we get 
\beq
\begin{gathered}
( L^t f, g )= (f, Lg) =  \int_{S\mcn} f(x, \theta) Lg(x, \theta) \sqrt{\det h} dx d\theta \\
= \int_{S\mcn} \int_\mbr f(x, \theta) g(s, \exp_x^{h}(s\theta)) \sqrt{\det h} ds dx d\theta 
\end{gathered}
\eeq
Let $y = \exp_x^{h}(s\theta)$. 
Because the Liouville measure is preserved by the geodesic flow, we have
\beq
\begin{gathered}
( L^t f, g )  = \int_{S\mcn} \int_\mbr f(\exp_{y}^{h}(-s\theta), \theta) g(s, y)\sqrt{\det h} ds dy d\theta 
\end{gathered}
\eeq
Thus 
\beq
L^tf(t, y) = \int_{S_y\mcn}  f(\exp_{y}^{h}(-t\theta), \theta)  d\theta
\eeq
Finally,
\beq
\begin{gathered}
L^tLf(t, y) 
=  \int_{S_y\mcn}\int_\mbr f(s,  \exp_{y}^{h}((s-t)\theta)) ds d\theta
\end{gathered}
\eeq
Let $\sigma = s - t$. We get 
\beqq\label{eq-inte1}
\begin{gathered}
L^tLf(t, y) 
=  \int_{S_y\mcn}\int_\mbr f(t +\sigma,  \exp_{y}^{h}(\sigma \theta)) d\sigma d\theta \\
= \int_{S_y\mcn}\int_0^\infty f(t +\sigma,  \exp_{y}^{h}(\sigma \theta)) d\sigma d\theta  + \int_{S_y\mcn}\int_0^\infty f(t -\sigma,  \exp_{y}^{h}(\sigma \theta)) d\sigma d\theta 
\end{gathered}
\eeqq
For the first integral, we use geodesic normal coordinate at $y$ so $x = \exp_y^{h}(\sigma\theta)$. Note that $\sigma = \dist^{h}(x, y)$. Using polar coordinates of the geodesic normal coordinate $w = (\sigma, \theta) \in T_x\mcn$, we have 
$\sigma^{n-1}d\sigma d\theta = \sqrt{\det h(y)} dw$. Then we change $w$ to $x$ via $x = \exp^h_y(w)$ and get  $\mcj(x, y) d\vol_h(x)$ where $\mcj(x, y)$ is the Jacobian of the change of variable. In particular, $\mcj(x, x) =  1.$ So we get  
\beqq
\begin{gathered}
 \int_{S_y\mcn}\int_0^\infty f(t +\sigma,  \exp_{y}^{h}(\sigma \theta)) d\sigma d\theta  
  = \int_{\mcn}  \frac{f(t + \dist^{h}(x, y), x)}{ (\dist^{h}(x, y))^{n-1}}\mcj(x, y) d\vol_h(x) \\
   = \int_\mbr  \int_\mcn \frac{\delta(t - t' - \dist^{h}(y, x))}{(\dist^{h}(y, x))^{n-1}} \mcj(x, y) f(t', x) d\vol_h(x) dt'
\end{gathered}
\eeqq
The second integral in \eqref{eq-inte1} can be treated similarly. 
\epf

\subsection{Globally hyperbolic spacetimes}
Consider globally hyperbolic spacetimes $(\mcm, g)$ with $g = -dt^2 + h(t, x, dx)$ as in \eqref{eq-ghmetric}.  Let $f(t, x, \tau, \xi) = \ha g(\tau, \xi)$ be the Hamilton function on $T^*\mcm.$ 
Let $\alpha(s) = (t(s), x(s), \tau(s), \xi(s))$ be a bicharacteristic curve. Then we get 
\beqq\label{eq-bisys1}
\begin{gathered}
\frac{dt}{ds}   = -\tau, \quad \frac{dx^j}{ds}  = \p_{\xi^j} h(t, x, \xi)\\
\frac{d\tau}{ds}   = -\p_{t}h(t, x, \xi), \quad \frac{d\xi_j}{ds}  = -\p_{x^j} h(t, x, \xi)
\end{gathered}
\eeqq
Comparing with \eqref{eq-bisys0}, we see that this system is not decoupled.  
Since we assumed that $(\mcm, g)$ is null-geodesically complete, we know that \eqref{eq-bisys1} has a unique solution defined for all $s\in \mbr.$  
It is convenient to use $t$ as the affine parameter. On $\Gamma^{lt}$, we see that $\tau \neq 0$ so we can write \eqref{eq-bisys1} as
\beqq\label{eq-bichar}
\begin{gathered}
\frac{ds}{dt} = -\frac{1}{\tau}, \quad \frac{dx^j}{dt} = -\frac{1}{\tau} \p_{\xi^j} h(t, x, \xi)\\
\frac{d\tau}{dt} =   \frac{1}{\tau}\p_{t}h(t, x, \xi), \quad \frac{d\xi_j}{dt} =  \frac{1}{\tau}\p_{x^j} h(t, x, \xi)
\end{gathered}
\eeqq
with initial condition
\beqq
s(t') = s_0, \quad x(t') = y, \quad \tau(t') = \tau_0, \quad \xi(t') = \eta.
\eeqq
The projection $(t, x(t)), t\in \mbr$ of $\gamma(t)$ to $\mcm$ is a light-like geodesic. Because $(\mcm, g)$ is globally hyperbolic, every hypersurface $\mcm_t \doteq \{t\}\times \mcn, t\in \mbr$ is a Cauchy surface. Thus the light-like geodesic $(t, x(t))$ intersect  with $\mcm_t$ at one point, see \cite{HaEl}. Consider the exponential map $\exp^g$  on future pointing light-like vectors 
\beqq\label{eq-texp}
 \exp^g_{(t',  y)}(\sigma(1, \theta))  = (t, x)
 \eeqq
 where $\theta = \eta^\#$ using the Riemmanian metric $h(t', x, dx)$ on $\mcn$ at $y$.  For $(t, x)\in \mcm$, we let $\mcl_\pm(t, x)$ be the future/past light cone at $(t, x)$
\beqq\label{def-ltcone}
\begin{gathered}
\mcl_\pm (t, x) = \{(t', x')\in \mcm: (t', x') = \exp^g_{(t, x)}s(\pm 1, \theta), \\
\text{  where $(\pm1, \theta)$  are light-like at $(t, x)$}, s \geq 0\}.
\end{gathered}
\eeqq 
We denote $\check{\mcl}_{\pm}(t, x) = \mcl_{\pm}(t, x)\backslash \{(t, x)\}$ the light cone minus the vertex. 
Since there is no conjugate point on $(\mcm, g),$  $\exp^g_{(t', y)}$  is a diffeomorphism from $T_y\mcn\backslash 0$ to $\check\mcl_+(t', y)$.  Using past pointing light-like vectors in \eqref{eq-texp}, we see that $\exp^g_{(t', y)}$ is also a  diffeomorphism from $T_y\mcn\backslash 0$ to $\check\mcl_-(t', y)$. Also, $\check\mcl_\pm(t', y)$ are smooth submanifolds of $\mcm.$ Below we mostly consider the exponential map to $\check\mcl_+(t', y)$ as the other case is similar. 
 
\begin{figure}[htbp]
\centering
\includegraphics[scale = 0.65]{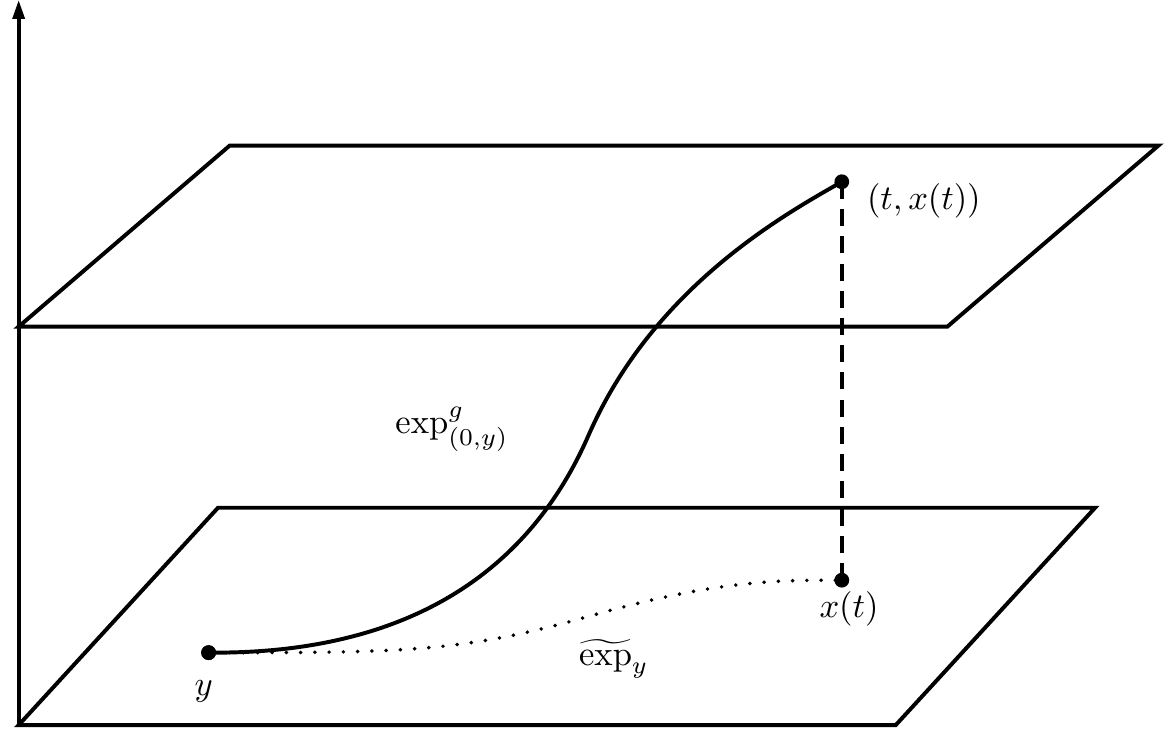}
\caption{Illustration of the map $\tilde{\exp}_{(t', y)}$ for $t' = 0.$}
\label{fig-exp}
\end{figure}

We need to be careful about the situation near the light cone vertex. In \eqref{eq-texp}, we denote $x = \tilde \exp_{(t', y)}(\sigma \theta).$ See Figure \ref{fig-exp}. For static spacetimes, the light cone can be identified with its projection to $\mcn$ and it suffices to look at $\tilde \exp_{(t', y)}$ which is $\exp^h$ in Section 4.1. In general, the projection does not seem to be always a diffeomorphism. 
 But for $x\in \mcn$ sufficiently close to $y$, this is true and we prove it next. It is convenient to compare it with the exponential map $\exp^{h_{t'}}$ on $(\mcn, h_{t'})$ with  $h_{t'} = h(t', x, dx)$. 

\begin{lemma}\label{lm-asymp}
Let $y\in \mcn$ and $x$ in a sufficiently small neighborhood $\mcv$ of $y$ such that there is no conjugate point on $(\mcv, h_{t}), t \in \mbr$. 
\begin{enumerate}
\item The map $\tilde \exp_{(t, y)}$ is a diffeomorphism from a small neighborhood of $0$ in $T_y\mcn$ to $\mcv$. For $x \in \mcv$, the function 
\beqq\label{eq-funcr}
\mcr(x, t, y) = |\tilde \exp_{(t, y)}^{-1}(x)|_{h_t}
\eeqq
is smooth in $x, t, y.$
\item  We have the following asymptotic expansion  
\beqq\label{eq-distasym}
\mcr(x, t, y) \sim  \sum_{k = 1}^\infty a_k(x, t, y)(\dist^{h_{t}}(x, y))^k 
\eeqq
where $a_k$ are smooth functions and $a_1 = 1$. Here, $\dist^{h_{t}}$ is the distance function on $(\mcv, h_{t})$.
\end{enumerate}
\end{lemma}
\bpf
Without loss of generality, we assume $t = 0$ and $\bar g = h(0, x, dx)$. 
We compare the map $\widetilde\exp$ with $\exp^{\bar g}$. Let $\bar\alpha(s) = (t(s),  x(s),  \tau(s),  \xi(s))$ be a bicharacteristics of $\bar f(t, x, \tau, \xi) = \ha(\tau^2 - \bar g(x, \xi))$. Then $\bar\alpha$ satisfies the corresponding system of \eqref{eq-bichar}
\beqq\label{eq-bichar1}
\begin{gathered}
\frac{d s}{dt} = -\frac{1}{\tau}, \quad \frac{d x^j}{dt} = -\frac{1}{\tau} \p_{\xi^j} \bar g(x, \xi)\\
\frac{d \tau}{dt} = 0, \quad \frac{d\xi_j}{dt} =  \frac{1}{\tau}\p_{x^j} \bar g(x, \xi)
\end{gathered}
\eeqq
with initial condition
\beqq
s(0) = s_0, \quad x(0) = y, \quad \tau(0) = \tau_0, \quad \xi(0) = \eta.
\eeqq
In this case, $\tau = \tau_0$ is a constant. The geodesic is denoted by $\bar \gamma(t) = x(t) = \exp^{\bar g}_y(t\eta)$. Here, we abused the notation of a covector $(\tau, \eta)\in T^*_{(s_0, y)}\mcm$ and its corresponding vector in $T_{(s_0, y)}\mcm.$ Now we consider the difference $\tilde \exp_y(t\eta) - \exp^{\bar g}_y(t\eta)$ on $S_y\mcn$. We write $\bar \alpha(t) = (\bar s(t), \bar x(t), \bar \tau(t), \bar \xi(t))$ and denote $\alpha(t) - \bar\alpha(t) = ((\delta s)(t), (\delta x)(t), (\delta \tau)(t), (\delta \xi)(t))$.  By using \eqref{eq-bichar} and \eqref{eq-bichar1}, they satisfy for $t$ small 
\beqq\label{eq-bichard}
\begin{gathered}
\frac{d (\delta s)}{dt} = -\frac{1}{\tau_0 + \delta \tau}, \quad \frac{d x^j}{dt} = -\frac{1}{\tau_0 + \delta \tau}t f_1\\
\frac{d (\delta \tau)}{dt} =  \frac{1}{\tau_0 + \delta \tau}t f_2, \quad \frac{d\xi_j}{dt} =  \frac{1}{\tau_0 + \delta \tau}t f_3
\end{gathered}
\eeqq
with zero initial condition, where $f_1, f_2, f_3$ are smooth functions of $t, \delta x, \delta \xi.$ By the small time existence result of nonlinear first order ODEs, we conclude that for $t \geq 0$ sufficiently small, $x(t) - \bar x(t) = F(t)$ as a map $T_y \mcn \rightarrow \mcn$ is  smooth. Also, $|F(t)|_{\bar g}\leq Ct$. We know that $\bar x(t)$ is invertible for $t$ small and $t = |\exp^{\bar g}_{y}\bar x|_{\bar g}$. Therefore, by using the inverse function theorem, we conclude that for $t$ sufficiently small, $x(t)$ is invertible and  
\beq
\tilde \exp_{(0, y)}^{-1} (x) = t (1 + \sum_{k = 1}^\infty a_k(y)t^k)
\eeq
where $t = \dist^{\bar g}(x, y).$ This gives the desired expansion. 
\epf

Now we consider the light ray transform in \eqref{eq-lorlray}. As discussed above, let $(x, \theta)\in S\mcn$, we have
\beqq\label{eq-lorlray1}
Lf(x, \theta) = \int_\mbr f(\exp^g_{(0, x)}(t(1, \theta)) dt.
\eeqq
On $\mcm$, we will use the measure $dg = d\vol_{h_t}(x)\wedge dt$ which in local coordinate reads $\sqrt{\det h(t, x)}dx\wedge dt$. On $S\mcn$, we will use the measure $d\vol_{h_0}(x) \wedge d\theta$ where $h_0 = h(0, x, dx)$. For a smooth submanifold $\Gamma$ of $\mcm$, we denote by $\delta_\Gamma$ the delta distribution supported on $\Gamma.$ 
\begin{prop}\label{prop-normal1}
Consider the light ray transform \eqref{eq-lorlray1} on globally hyperbolic manifolds of dimension $n+1, n\geq 2$ without conjugate points. Consider the Schwartz kernel of the normal operator $N = L^tL$.
\begin{enumerate}
\item Away from $t = t'$, the Schwartz kernel can be expressed as 
\beqq\label{eq-normal2}
\begin{gathered}
K_N(t, x, t', x') =    \frac{\delta_{\check \mcl_+(t, x)}(t', x')}{(t' - t)^{n-1}} \mcj(t, x, t', x')  \text{ for $t' > t$}\\
K_N(t, x, t', x') =      \frac{\delta_{\check \mcl_-(t, x)}(t', x')}{(t - t')^{n-1}} \mcj(t, x, t', x')   \text{ for $t' < t$}.
\end{gathered}
\eeqq
in which $(t, x), (t', x')\in \mcm$ and $\mcj(t, x, t', x')$ is a smooth non-vanishing function. 
\item For $(t', x') \in\mcm$ close to $(t, x)\in \mcm$, the Schwartz kernel can be expressed as 
 \beqq\label{eq-normal1}
\begin{gathered}
K_N(t, x, t', x') = \frac{\delta(t - t' - \mcr(x', t, x))}{(\mcr(x', t, x))^{n-1}}\mcj(t, x, x')  \\
+  \frac{\delta(t - t' + \mcr(x', t, x))}{(\mcr(x', t, x))^{n-1}} \mcj(t, x, x')
\end{gathered}
\eeqq
 where $\mcr$ is defined in \eqref{eq-funcr} and $\mcj(t, x, x')$ is a smooth non-vanishing function in $(t, x)$ and $x'$. 
\end{enumerate}
\end{prop}
\bpf
We proceed as in Proposition \ref{prop-normal0}.
First, we have
\beq
\begin{gathered}
( L^t f, g )= \int_{S\mcn} f(x, \theta) Lg(x, \theta) d\vol_{h_0}(x) d\theta 
= \int_{S\mcn} \int_\mbr f(x, \theta) g(\exp^g_{(0, x)}(s(1, \theta)) ds d\vol_{h_0}(x) d\theta 
\end{gathered}
\eeq
We set $(s, y)  = \exp^g_{(0, x)}(s(1, \theta))$. Because the exponential map is a diffeomorphism, we can write $(0, x) = \exp^g_{(s, y)}(-s(1, \tilde\theta))$ for some $\tilde \theta \in T_y\mcn$, $h(s, \tilde \theta) = 1$. We denote $x  = \Pi(-s, y, \tilde \theta) = \tilde \exp^g_{(s, y)}(-s(1, \tilde\theta)), \theta  = \Theta(-s, y, \tilde \theta) =  \frac{d}{dt}\tilde \exp^g_{(s, y)}(-t(1, \tilde\theta))|_{t = s}$. 
 We thus get 
\beq
\begin{gathered}
( L^t f, g )  = \int_{S\mcn} \int_\mbr f(\Pi(-s, y, \tilde \theta), \Theta(-s, y, \tilde \theta)) \tilde \mcj(s, y, \tilde \theta) g(s, y) ds d\vol_{h_0}(y) d\tilde \theta 
\end{gathered}
\eeq
where $\tilde J(s, y, \tilde \theta)$ is a smooth function induced by the Jacobian of the change of variable. Thus 
\beq
L^tf(t, y) = \int_{S_y\mcn}   f(\Pi(-t, y,  \theta), \Theta(-t, y,  \theta)) \tilde \mcj(t, y,  \theta) d \theta 
\eeq
Next, we have 
\beq
\begin{gathered}
L^tLf(t, y)  
=  \int_{S_y\mcn}\int_\mbr f(\exp^g_{(0, \Pi(-t, y, \theta))}(s(1, \Theta(-t, y, \theta))))\tilde \mcj(t, y, \theta) ds d\theta\\
 = \int_{S_y\mcn}\int_\mbr f(\exp^g_{(t, y)}((s - t)(1, \theta))\tilde \mcj(t, y, \theta) ds d\theta
\end{gathered}
\eeq
Let $\sigma = s - t$. We get 
\beqq\label{eq-llt}
\begin{gathered}
L^tLf(t, y) 
=  \int_{S_y\mcn}\int_\mbr f(\exp_{(t, y)}^g(\sigma(1, \theta))) \tilde \mcj(t, y, \theta) d\sigma d\theta \\
=  \int_{S_y\mcn}\int_0^\infty f(\exp_{(t, y)}^g(\sigma(1, \theta))) \tilde \mcj(t, y, \theta) d\sigma d\theta +  \int_{S_y\mcn}\int_0^\infty f(\exp_{(t, y)}^g(-\sigma(1, \theta))) \tilde \mcj(t, y, \theta) d\sigma d\theta
\end{gathered}
\eeqq

Consider part (1). For the first integral in \eqref{eq-llt}, we let $(t', x) = \exp_{(t, y)}^g(\sigma(1, \theta))$ so $t' = t + \sigma$. 
Since $\exp^g_{(t, y)}: T_y\mcn\backslash\{0\} \rightarrow \check\mcl_+(t, y)$ is a diffeomorphism, we again use polar coordinate on $T_y\mcn$ and the exponential map to get $\sigma^{n-1} d\sigma d\theta = \mci(t', x, t, y) dt' d\vol_{h_{t'}}(x)|_{\check\mcl_+(t, y)}$ where the right hand side denotes the restriction of the measure on $\check\mcl_+(t, y)$. Let $\mci$ be the Jacobian factor and note that $\theta$ is determined by $(t', x)$ and $(t, y)$. 
We get  
\beq
\begin{gathered}
  \int_{S_y\mcn}\int_0^\infty f(\exp_{(t, y)}^g(\sigma(1, \theta))) \tilde \mcj(t, y, \theta) d\sigma d\theta  
 =  \int_{\mcl_+(t, y)} \frac{1}{(t' - t)^{n-1}} \mcj(t, y, t', x) f(t', x)dt' d\vol_{h_{t'}}(x)  
\end{gathered}
\eeq
where $\mcj = \tilde \mcj \cdot \mci.$ 
The second integral in \eqref{eq-llt} can be treated similarly. This completes the proof of part (1).

For part (2), we can use Lemma \ref{lm-asymp} and write
\beq
\begin{gathered}
  \int_{S_y\mcn}\int_0^\infty f(\exp_{(t, y)}^g(\sigma(1, \theta))) \tilde \mcj(t, y, \theta) d\sigma d\theta   =   \int_{S_y\mcn}\int_0^\infty f(t' + \sigma, \tilde \exp_{(t, y)}(\sigma \theta)) \tilde \mcj(t, y, \theta) d\sigma d\theta
\end{gathered}
\eeq
Since $\tilde \exp_{(t, y)}: T_y \mcn\rightarrow \mcn$ is locally a diffeomorphism. We can use polar coordinate for $T_y\mcn$ to get  $\sigma^{n-1} d\sigma d\theta = \mci(t, y, x) d\vol_{h_{t'}}(x)$ where $\mci$ be the Jacobian factor. In this case, $\sigma = \mcr(x, t, y)$ and $t'$ is determined by $t, x, x'$. So we get  
\beq
\begin{gathered}
  \int_{S_y\mcn}\int_0^\infty f(\exp_{(t, y)}^g(\sigma(1, \theta))) \tilde \mcj(t, y, \theta) d\sigma d\theta  
   = \int_{\mcn} \frac{f(t + \mcr(x, t, y), x)}{\mcr(x,t, y)^{n-1}}   \mcj(t, y, x) d\vol_{h_{t'}}(x)
\end{gathered}
\eeq
The second integral in \eqref{eq-llt} is similar. This completes the proof of the proposition. 
\epf

\section{Proof of Theorem \ref{thm-main}}\label{sec-mainres}
As indicated by Proposition \ref{prop-normal1}, we define two Lagrangian submanifolds $\La_0, \La_1$ as follows. Let  $f = \ha g(\zeta, \zeta)$ and $\Sigma = \{(z, \zeta, z, -\zeta) \in T^*(\mcm\times\mcm)\backslash 0:  f(z, \zeta) = 0\}$. Then we let 
\beqq\label{eq-laggen}
\begin{gathered}
\text{$\La_0 = N^*\diag$   the conormal bundle of the diagonal of $\mcm\times \mcm$}, \\
\text{$\La_1 = $ the flow out of $\Sigma$ under $H_f$}  
\end{gathered}
\eeqq 
Observe that $\La_0$ intersect $\La_1$ cleanly at $\Sigma$ which has co-dimension one in $T^*(\mcm\times \mcm)$. Also, away from $\La_0$, $\La_1$ is the conormal bundle of the set $\{(t, x, t', x')\in \mcm\times \mcm: (t', x')\in \check\mcl_\pm(t, x)\}$.   

\bpf[Proof of Theorem \ref{thm-main}]
Away from $\La_0$, one can see from Proposition \ref{prop-normal1} part (1) that $N \in I^{-n/2}(\mcm; \La_1)$. Now we examine the kernel near $\La_0\cap \La_1$ by using Proposition \ref{prop-normal1} part (2). 
We use the kernel expression in a local coordinate $(t, x, t', x')$ for $\mcm\times \mcm.$ By introducing a smooth partition of unity, it suffices to assume that the Schwartz kernel is compactly supported in this coordinate patch. Let $\chi(s)$ be a smooth cut-off function which is supported in a sufficiently small neighborhood of $0$. We consider the kernel
\beqq\label{eq-kn}
\begin{gathered}
K_N(t, x, t', x') \chi(t - t')\\
=\frac{\delta(t - t' - \mcr(x', t, x))}{(\mcr(x', t, x))^{n-1}}\mcj(t, x, x')  \chi(t - t')
+  \frac{\delta(t - t' + \mcr(x', t, x))}{(\mcr(x', t, x))^{n-1}} \mcj(t, x, x') \chi(t - t')\\
 = \int_\mbr  e^{i(t - t' - \mcr(x', t, x))\tau} \frac{1}{(\mcr(x', t, x))^{n-1}} \mcj(t, x, x') \chi(t - t') d\tau  \\
  +\int_\mbr   e^{i(t - t' + \mcr(x', t, x))\tau} \frac{1}{(\mcr(x', t, x))^{n-1}} \mcj(t, x, x') \chi(t - t')  d\tau   
 = I_- + I_+  
\end{gathered}
\eeqq 
To simplify the calculation, we will  use geodesic normal coordinate based at $(t', x')$ so $g$ agrees with the Minkowski metric at $(t', x').$ Then we can write  $(t-t', x-x') = \sigma(1, \theta), \theta \in S_{x'}\mcn$ with $\sigma = t - t'$ for $t\geq t'.$ For $t\leq t'$, we can write $(t - t', x - x') = \sigma(-1, \theta)$ with $\sigma = t' - t \geq 0, \theta \in S_{x'}\mcn$. Note that $\mcr(x', t, x)$ is only defined for $x$ close to $x'$. But because of the $\chi(t'-t)$ factor, we can extend $\mcr(x', t, x)$ smoothly outside the neighborhood of $x'$.  So we assume that $\mcr(x', t, x) = \sigma$ for all $(t, x)\in \mbr^{n+1}$ in geodesic normal coordinate  which does not affect $K_N\chi$.

For $I_-$, we write
\beq
\begin{gathered}
I_- = \int_{\mbr^{n+1}} e^{i(t - t')\tau}  e^{i(x - x')\xi} A_-(x', t, \tau, \xi)  \mcj(t, x, x') \chi(t - t') d\tau d\xi  
\end{gathered}
\eeq
 where 
 \beq
 \begin{gathered}
 A_-(x', t, \tau, \xi) = \int_{\mbr^{n}} e^{i(x' - y)\xi}  e^{-i\mcr(x', t, y)\tau}   \frac{1}{(\mcr( x', t, y))^{n-1}}  dy\\
 = \int_{\mbr^{n}} e^{i z \xi}  e^{-i\mcr(x', t, z)\tau}   \frac{1}{(\mcr( x', t, z))^{n-1}} dz
 \end{gathered}
 \eeq
where we changed variables $z = y - x'$ and write $\mcr(x', t, z) = \mcr(x', t, z + x')$ in the second line. We  write  $(\sigma, z) = \sigma(1, \theta), \theta\in S_{x'}\mcn$ and $\sigma = t - t'$ for $t\geq t'$. We get
  \beqq\label{eq-aminus}
 \begin{gathered}
 A_-(x', t, \tau, \xi)  
  = \int_{0}^\infty \int_{\mbs^{n-1}} e^{-i  \sigma \theta \cdot  \xi }  e^{-i\sigma \tau} \sigma^{-(n-1)} \sigma^{n-1} d\sigma d\theta  
 \end{gathered}
 \eeqq
 For $t\leq t'$, we can write $(\sigma, z) = \sigma(-1, \theta)$ with $\sigma = t' - t \geq 0$. So  
  \beq 
 \begin{gathered}
 A_-(x', t, \tau, \xi)  
  = \int_{0}^\infty \int_{\mbs^{n-1}} e^{-i  \sigma \theta \cdot  \xi }  e^{-i\sigma \tau} \sigma^{-(n-1)} \sigma^{n-1} d\sigma d\theta  
 \end{gathered}
 \eeq 
 Thus the expression of  $A_-$ in \eqref{eq-aminus} holds for all $t, t'$. For $I_+,$ we repeat the argument and get  
 \beq
\begin{gathered}
I_+ = \int_{\mbr^{n+1}} e^{i(t - t')\tau}  e^{i(x - x')\xi} A_+(x', t, \tau, \xi)  \mcj(t, x, x')\chi(t - t')  d\tau d\xi  
\end{gathered}
\eeq
 where 
 \beq
 \begin{gathered}
 A_+(x',t, \tau, \xi) = \int_{\mbr^{n}} e^{i(x' - y)\xi}  e^{i\mcr(x', t, y)\tau}   \frac{1}{(\mcr(  x', t, y))^{n-1}} dy\\ 
  = \int_{0}^\infty \int_{\mbs^{n-1}} e^{-i\sigma \theta \cdot \xi}  e^{i\sigma \tau} \sigma^{-(n-1)}   \sigma^{n-1}d\sigma d\theta  
 \end{gathered}
 \eeq 
Therefore,   
 \beqq\label{eq-apm} 
 A_-(x',t, \tau, \xi)+  A_+(x',t, \tau, \xi) =  \int_{0}^\infty \int_{\mbs^{n-1}} e^{-i\sigma \theta \cdot \xi} ( e^{i\sigma \tau} + e^{-i\sigma \tau})   d\sigma d\theta
 \eeqq 
Now we can use the calculation in \cite{LOSU2}  to conclude that 
  \beqq\label{eq-apm2} 
 \begin{gathered}
 A_-(x',t, \tau, \xi)+  A_+(x',t, \tau, \xi)   = C_n    \frac{(|\xi|^2 - \tau^2)_+^{\frac{n-3}{2}}}{|\xi|^{n-2}}   
  \end{gathered}
 \eeqq
 Thus 
 \beq 
\begin{gathered}
K_N(t, x, t', x') \chi(t - t') = \int_{\mbr^{n+1}} e^{i(t - t')\tau}  e^{i(x - x')\xi} C_n   \frac{(|\xi|^2 - \tau^2)_+^{\frac{n-3}{2}}}{|\xi|^{n-2}}    \mcj(t, x, x')\chi(t - t')  d\tau d\xi  
\end{gathered}
\eeq  
 
 We can argue as in Section \ref{sec-pair} to see that $K_N$ in \eqref{eq-kn} is a paired Lagrangian distribution associated with $\La_0, \La_1$. The principal symbols on $\La_0, \La_1$ can be computed as in Section \ref{sec-pair} as well. The difference is the $\mcj(t, x, x')\chi(t - t')$ factor so the symbols are non-vanishing. 
 \epf

\bpf[Proof of Theorem \ref{thm-mainest}]
The proof follows from Theorem \ref{thm-main} and the Sobolev estimate Theorem \ref{thm-sobolev}. 
\epf

We remark that one can find a parametrix using the ellipticity of the paired Lagrangian distribution as in Proposition \ref{prop-para1} but with a $I^{-1/2}(\mcm; \La_1)$ remainder. We will not carry out the details which can be found in \cite{Wan} and \cite{GrUh0}. 
 
\section{Light-like singularities}\label{sec-det}
We are interested in the determination of light-like singularities using the normal operator.  Suppose $\Gamma \subseteq T^*\mcm\backslash 0$ be a conic Lagrangian submanifold and $\Gamma \subseteq \Gamma^{lt}$. We assume that $f\in I^\mu(\mcm; \Gamma), \mu\in \mbr$ is a compactly supported Lagrangian distribution with non-trivial principal symbol. Since $N \in I^{-n/2, n/2 -1}(\mcm\times \mcm; \La_0, \La_1)$  and $f\in I^\mu(\mcm; \Gamma)$, we can use the composition result Proposition 2.3 of  \cite{GrUh1} to conclude that 
\beq
Nf \in I^{\mu -n/2 + 1/2}(\mcm; \Gamma_1)
\eeq
where $\Gamma_1$ is the flow out of $\Gamma$ under $H_f$. 
We will need the principal symbol, which can be seen by examining the arguments of Proposition 2.1-2.3 in \cite{GrUh1}.   As in Proposition 2.1 of \cite{GrUh1}, we can assume that $\mcm = \mbr^{n+1}$ with coordinates $x = (x_1, x'), x_1\in \mbr, x'\in \mbr^n, \Gamma_0 = T_0^*\mbr^{n+1}\backslash 0$ and $\Gamma^{lt} = \{(x, \xi): \xi_1 = 0\}$. Thus 
$
\Gamma = \{(x_1, 0; 0, \xi'): x_1\in \mbr, \xi'\in \mbr^{n}\backslash 0\}. 
$
In this coordinate, we can write
\beqq\label{eq-f}
f(x) = \int_{\mbr^n} e^{ix'\xi'} a(x_1, \xi') d\xi', \quad a\in S^{\mu - (n+1)/4}(\mbr \times (\mbr^{n}\backslash 0))
\eeqq
and  
\beqq\label{eq-Nf1}
Nf(x) = \int_{\mbr^{n+1}\times \mbr^{n+1}} e^{i(x - y)\cdot \theta} b(x, y; \theta', \theta_1)f(y)d\theta dy
\eeqq
where $b\in S^{p+1/2, l-1/2}(\mbr^{n+1}\times \mbr^{n+1}\times (\mbr^{n}\backslash 0)\times \mbr)$. 
Applying \eqref{eq-f} to \eqref{eq-Nf1}, we obtain 
\beq 
Nf(x) =  \int_{\mbr^n} \int_{\mbr^{n+1}\times \mbr^{n+1}} e^{i(x - y)\cdot \theta} e^{iy'\xi'}a(y_1, \xi')b(x, y; \theta', \theta_1) d\theta dy d\xi'
\eeq 
Applying  the stationary phase in $y', \theta'$ and relabelling  $\theta_1$ by $\xi_1$, we get 
\beq 
\begin{gathered}
Nf(x) = \int_{\mbr \times \mbr^n} e^{ix \xi} (\int_\mbr e^{-iy_1\xi_1}a(y_1, \xi') b(x, (y_1, 0); \xi', \xi_1)   dy_1) d\xi_1d\xi'\\ 
= \int_{\mbr^n} e^{ix' \xi'} [ \int_\mbr e^{i(x_1 - y_1)\xi_1}a(y_1, \xi') b(x, (y_1, 0); \xi', \xi_1)   dy_1  d\xi_1] d\xi'
\end{gathered}
\eeq 
We see that 
this is an element of $I^{p + r}(\mcm; \Gamma_1)$ when the integral in $y_1, \xi_1$ is well-defined. In fact, the leading order terms of the integral in the square bracket as $|\xi'|\rightarrow \infty$ is the product of the principal symbols $\sigma(f)$ and $\sigma_{\La_1\backslash \La_0}(N)$ (when $y_1\neq x_1$). 
In case of the normal operator $N = L^tL$, we see from \eqref{eq-kker}, \eqref{eq-kn3} and the calculations in Section \ref{sec-mainres} that the symbol $\sigma_{\La_1\backslash \La_0}(N)$ is singular at $\La_0\cap \La_1$. However, the integral in the square bracket exists for $n=2$ and exists in the principal value sense when $n = 3$. Thus, we have
\begin{prop}\label{prop-composition}
Let $f\in I^{\mu}(\mcm; \Gamma)$. Consider the normal operator in Theorem \ref{thm-main1} with $n = 2, 3$. Then $Nf \in I^{\mu -1}(\mcm; \Gamma_1)$ and the principal symbol is given by
\beq
\sigma(Nf)(z, \zeta) = \int \sigma_{\La_1\backslash \La_0}(N)(z, \zeta, w, \iota)\sigma(f)(w, \iota)
\eeq
where the integration is over the bicharacteristics from $(z, \zeta)$ and exists in the sense of principal value for $n=3.$
\end{prop}
 
\bpf[Proof of Theorem \ref{thm-main1}]
We start with a compactly supported distribution  $f\in I^\mu(\mcm; \Gamma)$ with principal symbol $a = \sigma(f) \neq 0$. 
We assume that $\supp(f)$ is sufficiently small and is contained in an open set $\mcu$ with $\overline\mcu$ compact.  Let  $\gamma_p(s), s\in \mbr$ be the bicharacteristics from $p \in T^*\mcu $  and $\gamma_p(0) = p.$ For any $p'= \gamma_p(s'), s\in \mbr$, we consider
\beq
\alpha(\gamma_p(s')) = \int_{\mbr} \sigma(N)(\gamma_p(s'), \gamma_p(s))a(\gamma_p(s))ds 
\eeq 
Note that $\alpha$ is constant for a fixed bicharacteristic. In fact $\alpha(\gamma_p(s'))$ is the principal symbol of $Nf$ at $\gamma_p(s')$. 
Now we choose a compactly supported distribution $h\in I^\mu(\mcm; \Gamma)$ with principal symbol $b(p) = a(\gamma_p(T)), p\in T^*\mcm$ with  $T >0$ fixed. Here, we see that $h$ is compactly supported and we can shrink the support of $f$ and choose $T$   sufficiently large so that $\pi(\supp(b)) \cap \pi( \supp(a)) = \emptyset$ where $\pi: T^*\mcm\rightarrow \mcm$ is the natural projection. Then we consider   
\beq
\begin{gathered}
\beta(\gamma_p(s'))= \int_\mbr \sigma(N)(\gamma_p(s'), \gamma_p(s)) b(\gamma_p(s))ds   
\end{gathered}
\eeq
which is also a constant for fixed bicharacteristic. Now we consider a distribution $f_0\in I^\mu(\mcm; \Gamma)$ with principal symbol 
\beq
\sigma(f_0)(\gamma_p(s)) = \beta(\gamma_p(s))a(\gamma_p(s)) - \alpha(\gamma_p(s))b(\gamma_p(s)) 
\eeq
Then according to Proposition \ref{prop-composition}, we get $\sigma(Nf_0)(\gamma_p(s')) =  0$ 
 in $I^{\mu - n/2 -1/2}(\mcm; \Gamma_1)$  so $Nf_0 \in I^{\mu -n/2 - 1/2}(\mcm; \Gamma_1), n= 2, 3$. This argument can be continued to give $f_j\in I^{\mu - n/2 +1/2 - j}(\mcm; \Gamma),$ $j\geq 1$ such that 
\beq
N(\sum_{j = 0}^J f_j ) \in I^{\mu -n/2 + 1/2 - J}(\mcm; \Gamma_1)
\eeq
Thus taking asymptotic summation $\tilde f \sim \sum_{j = 0}^J f_j$, we constructed a distribution $\tilde f \in I^{\mu - n/2 - 1/2}(\mcm; \Gamma), n = 2, 3$ such that $\WF \tilde f \neq \emptyset$ but $N \tilde f \in C^\infty(\mcm).$ This completes the proof of the theorem. 
\epf

The above result suggests that using singularities of the normal operator to determine light-like singularities, one needs to consider more restrictive Lagrangian distributions.  We give one example.  
\begin{theorem}\label{thm-main2}
Suppose $f\in I^\mu(\mcm; \Gamma)$ with compact support and $n=2, 3$. Suppose that the principal symbol of $f$ is non-positive or non-negative along each bicharacteristics in $T^*\mcm$. If $Lf \in C^\infty(\mcc)$, then $f\in I^{\mu-1}(\mcm; \Gamma).$
\end{theorem} 
\bpf
This is immediate from Proposition \ref{prop-composition} and the kernel representation \eqref{eq-normal2} in Proposition \ref{prop-normal1}.  
 \epf

\section{Timelike caustics on static spacetimes}\label{sec-conjugate}
In this section, we analyze the effects of conjugate points in the setting of standard static spacetimes when the exponential map can be reduced to that of a Riemmanian manifold. The conjugate points in geodesic ray transforms in the Riemmanian setting was studied in Stefanov and Uhlmann \cite{StUh} where fold type conjugate points are analyzed. Later in   \cite{HoUh}, the regular conjugate points are analyzed using the clean FIO calculus. Here, we will use the results in \cite{StUh}. For Lorenzian manifolds, Rosquist proved analogues of Warner's result \cite{War} about classification of conjugate points, see Theorem 3.8 and 3.9 of \cite{Ros} for the statement of time-like conjugate points. See also \cite{Sze} for the generalization to pseudo-Riemannian manifolds. We will consider the regular conjugate points of fold type and we recall the preliminaries below, see \cite{StUh, Ros, War} for details.

Let $(\mcx, h)$ be a Riemannian or Lorentzian manifold of dimension $n$. Let $p\in \mcx$ and $\exp_p: T_p\mcx\rightarrow \mcx$ be the exponential map on $(\mcx, h)$. The tangent conjugate locus $S(p)$ of $p$ is the set of all $v\in T_p\mcx$ such that $d_v\exp_p v$ is not an isomorphism.  Such vectors $v$ are called conjugate vectors at $p$. In the Lorentzian case, one can classify the conjugate vectors to time-like, light-like and space-like. The kernel of $d\exp_p$ is denoted by $N_p(v) \subseteq T_v(T_p \mcx)$ and can be identified as a set of $T_p \mcx$. The image of $S(p)$ under $\exp_p$ is called conjugate locus $\Sigma(p)$. For $v\in S(p), p\in \mcm$, let $q = \exp_{p}(v)$. We denote by $\Sigma$ the set of all conjugate points pair $\Sigma = \{(p, q): q\in \Sigma(p)\}$. Similarly, $S = \{(p, v): v\in S(p)\}$. A regular conjugate vector $v$ is that there is a neighborhood $\mcb$ of $v$ such that any radial ray in $\mcb$ contains at most one conjugate point. The regular conjugate locus is a dense open subset of $\Sigma(p)$ and $n-1$ dimensional manifold. Among regular conjugate points, the fold type conjugate point is defined as follows. 
\begin{definition}\label{def-fold}
A regular conjugate vector $v_0$ at $p_0$ is called fold type if $N_{p_0}(v_0)$ is $1$-dimensional and transversal to $S(p_0)$.
\end{definition}
In this case, there exists local coordinate $\xi$ near $v_0$, $y$ near $q_0$ such that $\exp_{p_0}$ is expressed by
\beq
y' = \xi', \quad y^n = (\xi^n)^2
\eeq
and that $S(p_0) = \{\xi^n = 0\}, N_{p_0}(v_0) = \text{span}(\p/\p \xi^n)$ and $\Sigma(p_0) = \{y^n = 0\}$. The fold condition is stable under small $C^\infty$ metric perturbations and small perturbations of $p_0.$

Now we consider standard static spacetimes $(\mcm, g)$ as in \eqref{eq-staticmetric}. Let $v_0$ be a regular time-like conjugate vector at $p_0\in \mcm$. Let $\exp^g$ be the exponential map on $(\mcm, g).$ As in \cite{StUh}, we will prove a local result and assume that the support of $f$ is such that $v_0$ is the only conjugate vector at $p_0$ such that $\exp_{p_0}^g (v) \in \supp f.$ Using \eqref{eq-staticmetric}, we write $p = (t, x), t\in \mbr, x\in \mcn$ and $v = (\beta, \theta), \beta\in \mbr, \theta\in T_x\mcn$. We consider the tangent conjugate locus $S(p)$ near $(p_0, v_0) = (t_0, x_0, \beta_0, \theta_0)$ and the corresponding conjugate locus   
 \beq
 \Sigma = \{(t, x, s, y)\in \mcm\times \mcm : (s, y) \in \Sigma(t, x), \text{$(s, y)$ close to $(t_0, x_0)$}\}.
 \eeq 
Locally, this is a $2n+1$ dimensional submanifold of $\mcm\times \mcm.$ 
  It is convenient to think in terms of the conjugate points on $(\mcn, h)$. Because the exponential map can be decomposed as in \eqref{eq-geo0}, we find that for $v = (\beta, \theta) \in T_{(t, x)}\mcm$, 
\beq
\exp^g_{(t, x)}(v) = (t + \beta, \exp^h_x\theta) 
\eeq
Therefore, $d_v\exp^g_{(t, x)}(v) = (1, d_\theta \exp^h_x\theta)$. If $v_0 = (\beta_0, \theta_0)$ is a regular conjugate vector of fold type at $(t_0, x_0)$, then $\theta_0$ is a regular conjugate vector of fold type for $(\mcn, h)$ at $x_0$.  Let $\Sigma^h = \{(x, y)\in \mcn\times \mcn: x \in \Sigma^h(y), \text{$y$ is close to $x_0$}\}$. We know from \cite[Theorem 2.1]{StUh} that 
  \beqq\label{eq-locus2}
 \begin{gathered}
 N^*\Sigma^h = \{(x, y, \xi, \eta): (x, y)\in \Sigma^h,  \xi = \sum_{i = 1}^n \eta_i \p\exp^{h, i}_y(\theta)/\p y, \\
 \eta \in \text{Coker}d_\theta \exp_y^h \theta, 
 \det d_\theta \exp^h_y(\theta) = 0\}. 
 \end{gathered}
 \eeqq
Thus, 
 \beqq\label{eq-locus1}
 \begin{gathered}
 N^*\Sigma = \{(t, x, s, y, \tau, \xi, \iota, \eta): (t, x, s, y)\in \Sigma, \tau= -\iota, 
 \xi = \sum_{i = 1}^n \eta_i \p\exp^{h, i}_y(\theta)/\p y, \\
 \eta \in \text{Coker}d_\theta \exp_y^h \theta, 
 \det d_\theta \exp^h_y(\theta) = 0\}. 
 \end{gathered}
 \eeqq

We are ready to state the main result of this section. 
\begin{theorem}\label{thm-main3}
Consider the light ray transform on standard static spacetime $(\mcm, g)$ as in \eqref{eq-staticmetric} with $n\geq 2$. Let $v_0$ be a regular time-like conjugate point at $p_0$ of fold type.  Let $\mcu\subset \mcm$ be an open set such that $v_0$ is the only conjugate vector $v$ at $p_0$ such that $\exp^g_{p_0}(v)\in \mcu$. Then the normal operator $N$ in Proposition \ref{prop-normal0} regarded as an operator $C_0^\infty(\mcu)\rightarrow C^\infty (\mcm)$ can be written as  $N = A + B$ 
where the Schwartz kernel $K_A\in  I^{-n/2, n/2- 1}(\mcm\times \mcm; \La_0, \La_1)$ with $\La_0, \La_1$ in \eqref{eq-laggen}, and $K_B\in I^{-n/2}(\mcm\times \mcm; N^*\Sigma)$ is a Lagrangian distribution associated with the conjugate pairs defined in \eqref{eq-locus1}. 
\end{theorem}
 \bpf
We use the calculation in Proposition \ref{prop-normal0}. In particular, from \eqref{eq-inte1}, we have 
\beq
\begin{gathered}
Nf(t, x) = L^tLf(t, x)  
 =  \int_{S_x\mcn}\int_\mbr f(t +\sigma,  \exp_{x}^{h}(\sigma \Theta)) d\sigma d\Theta \\
\end{gathered}
\eeq
We want to change variables $y = \exp_x(\sigma \Theta)$. For $(x, y)$ away from $\Sigma^h$, the exponential map is an isomorphism. Let $K_A$ be the Schwartz kernel of $N$ away from $\Sigma^h$. Then the calculation in Proposition \ref{prop-normal0} works and Proposition \ref{prop-normal1} shows  that $K_A$ is a paired Lagrangian distribution. 

Next we consider the part of the kernel near $\Sigma^h$, denoted by $K_B$. For this, we can use the calculations in \cite[Section 6]{StUh}. Since the calculations are quite involved, we will only sketch the main ingredients below and refer the reader to \cite{StUh} for details.  Consider $y  = \exp_x^h(\theta), \theta = \sigma \Theta\in T_x \mcn$ for $(x, \theta)$ close to $(x_0, \theta_0)$. Because $\theta_0$ is of fold type at $x_0$, the exponential map is $2$-to-$1$ near $S(x)$ for $x$ close to $x_0$. Thus, one should analyze on each sides of $S(x)$ so $\exp^h$ is $1$-to-$1$ as in \cite[Section 6]{StUh}. Let $y = (y', y^n)$ be semi-geodesic coordinate near $\Sigma^h(x_0)$ such that $y_0 = \exp_{x_0}^h(\theta_0)$ and $\Sigma^h(x)$ is given by $y^n = 0$. We consider the side corresponding to $y_n>0.$ Then let $\theta = (\theta', \theta^n)$ be boundary normal coordinates to $S(x)$ so $S(x)$ is given by $\theta^n = 0$ and we consider the side corresponding to $\theta^n>0$. In such coordinates, 
\beq
y^n = (\theta^n)^2 D_1 (1 + O(\theta^n))
\eeq
where $D_1$ and $O(\theta^n)$ denote smooth functions of $\theta, x$, and 
\beq
\det d_\theta \exp^h_x(\theta) = (y^n)^\ha( D_2 + O((y^n)^\ha))
\eeq
where $D_2,  O((x^n)^\ha)$ are also smooth of $y, x$. Using these results, 
 we have 
\beq
\begin{gathered}
B f(t, x)  =  \int_0^\infty\int_{\mbr^{n-1}}   |\theta|^{1-n}(\det d\exp^{h}_{x}(\theta))^{-1} f(t \pm |\theta|, y) dy' dy^n \\
 =  \int_0^\infty\int_{\mbr^{n-1}} \frac{D_3}{(y^n)^\ha} f(t \pm D_4, y)dy' dy^n
\end{gathered}
\eeq
where $D_3, D_4$ are smooth functions of $y, x.$ We further get 
\beq
\begin{gathered}
B f(t, x)   
=\int_\mbr  \int_0^\infty\int_{\mbr^{n-1}} \int_\mbr e^{i(t \pm D_4 - t') \tau}D_3 (y^n)^{-1/2}f(t', y)   d\tau dy' dy^n  dt'
\end{gathered}
\eeq
This shows that the kernel has a conormal singularity   at $y^n = 0$  and $t' = t \pm D_4(y, x)$ which corresponds to $\Sigma$. The analogue of $B$ corresponding to $y^n<0$ is similar. By taking Fourier transform in $y^n$ as in \cite[Section 7]{StUh}, we see that $B$ is an Fourier integral operator with a symbol of order $-1/2.$ So the order of the FIO is $-1/2 + 2/2 - (2n+2)/4 = -n/2$. Thus $K_B\in I^{-n/2}(\mcm\times \mcm; N^*\Sigma).$
 \epf
 
 One can use Theorem \ref{thm-main3} to formulate and prove analogous results as in Section \ref{sec-det}. We leave the details to the interested readers.

\section*{Acknowledgments}
The author thanks  Prof.\ Gunther Uhlmann and  Prof.\ Plamen Stefanov for  helpful comments.


\end{document}